\documentclass[11pt,a4paper]{article}
\usepackage[onehalfspacing]{setspace}
\usepackage[a4paper,margin=0.8in]{geometry}

\usepackage{graphicx}
\usepackage[english]{babel}
\usepackage{amsmath, amssymb, amsthm, amsfonts}
\usepackage{mathtools}
\usepackage{verbatim}
\usepackage{xurl}
\usepackage{multirow}
\usepackage{booktabs}
\usepackage{enumitem}

\usepackage{hyperref}
\usepackage{xcolor}
\usepackage[capitalise]{cleveref}
\usepackage{natbib}
\usepackage{authblk}
\usepackage{cases}
\usepackage{csquotes}

\usepackage{algorithm}
\usepackage{algpseudocode}
\algrenewcommand\algorithmicrequire{\textbf{Storage:}}
\algrenewcommand\algorithmicensure{\textbf{Input:}}

\renewcommand{\l}{\left}
\renewcommand{\r}{\right}
\newcommand{\E}{\mathbb{E}}
\newcommand{\iidsim}{\stackrel{\mathrm{iid}}{\sim}}

\renewcommand{\P}{\mathbb{P}}

\newcommand{\floor}[1]{\left\lfloor #1 \right\rfloor}
\newcommand{\ceil}[1]{\left\lceil #1 \right\rceil}

\newcommand{\rev}{}

\newtheorem{theorem}{Theorem}
\newtheorem{proposition}[theorem]{Proposition}

\newtheorem{corollary}[theorem]{Corollary}
\theoremstyle{definition}

\title{Online jump and kink detection in segmented linear regression: Statistical optimality meets computational efficiency}

\author[1]{Annika Hüselitz}
\author[1,2]{Housen Li}
\author[1,2]{Axel Munk}
\affil[1]{Institute for Mathematical Stochastics, University of G\"ottingen}
\affil[2]{Cluster of Excellence ``Multiscale Bioimaging: from Molecular Machines to Networks of Excitable Cells'' (MBExC), University of G\"ottingen}
\date{\today}

\begin{document}

\maketitle

\begin{abstract}
We consider the problem of sequential (online) estimation of a single change point in a piecewise linear regression model under a Gaussian setup. We demonstrate that certain  CUSUM-type statistics attain the minimax optimal rates for localizing the change point. Our minimax analysis unveils an interesting phase transition from a \emph{jump} (discontinuity in function values) to a \emph{kink} (a change in slope). Specifically, for a jump, the minimax rate is of order $\log (n) / n$ , whereas for a kink it scales as $\bigl(\log (n) / n\bigr)^{1/3}$, given that the sampling rate is of order $1/n$. {\rev We further  introduce an online algorithm based on these detectors, which optimally identifies both a jump and a kink, and is able to distinguish between them.} Notably, the algorithm operates with constant computational complexity and requires only constant memory per incoming sample. Finally, we evaluate the empirical performance of our method on both simulated and real-world data sets.  An implementation is available in the R package \texttt{FLOC} on GitHub.
\end{abstract}

Keywords: Sequential detection, minimax rate, efficient computation, Covid-19, two phase regression.

\medskip

\href{https://mathscinet.ams.org/mathscinet/msc/msc2020.html} {MSC2020 subject classifications:} Primary 62L12, 62C20; secondary 62J20. 

\section{Introduction}

Sequential change point detection, also known as online or quickest change point detection, is a classic topic in statistics with roots in ordnance testing and quality control \citep{shewhart1931economic, Wald45, Ansc46}. The primary goal is to monitor a (random) process and issue an alert when the distribution of this process deviates significantly from its historical pattern. Applications range widely, including social network monitoring \citep{Che19}, quality control in industrial processes \citep{amiri2012change}, cybersecurity \citep{TaPoSo12}, and seismic tremor detection \citep{li18high}, among others.

In early work,  \citet{Page54} introduced the celebrated CUSUM statistic based on the  the log-likelihood ratio between two known distributions (representing the control and the anomaly). \citet{Lord71} later established the asymptotic minimax optimality of CUSUM regarding detection delay under constraints on average run lengths. From a Bayesian perspective, the Shiryaev--Roberts procedure \citep{shiryaev1961problem,Shir63,Rob66} is similar to CUSUM and has several optimality properties (including the aforementioned minimax optimality as well as Bayesian optimality; see e.g.\ \citealp{Poll85,PoTa10}). Subsequent developments have expanded this scope to address unknown anomaly distributions and various models (see e.g.\ \citealp{Sie85}, \citealp{Lai01}, \citealp{Tar20} and \citealp{WaXi24} for an overview). 

Beyond statistical optimality, recent advancements emphasize computational and memory efficiency as essential \emph{design criteria} for sequential change point analysis, especially in modern data contexts (cf.\ arguments in \citealp[Section~1]{ChWaSa22}). This consideration becomes especially pertinent in scenarios where evaluating each candidate change point is computationally intensive \citep{KoBuLiMu23} and where memory resources are limited. Although computationally and memory-efficient methods are urgently needed for modern large-scale data analysis, they remain largely unexplored, with only a few exceptions \citep{kovacs2020optimistic,RoEcFeRi23, romano2023log, WaRoEcFe24}. A more detailed discussion of the literature from the perspective of this paper, with particular focus on minimax optimality and computational efficiency, will be presented in \Cref{ss:literature}. 

\subsection{Model and main results}\label{section:model}
\begin{figure}
    \centering    \includegraphics[width=.65\linewidth]{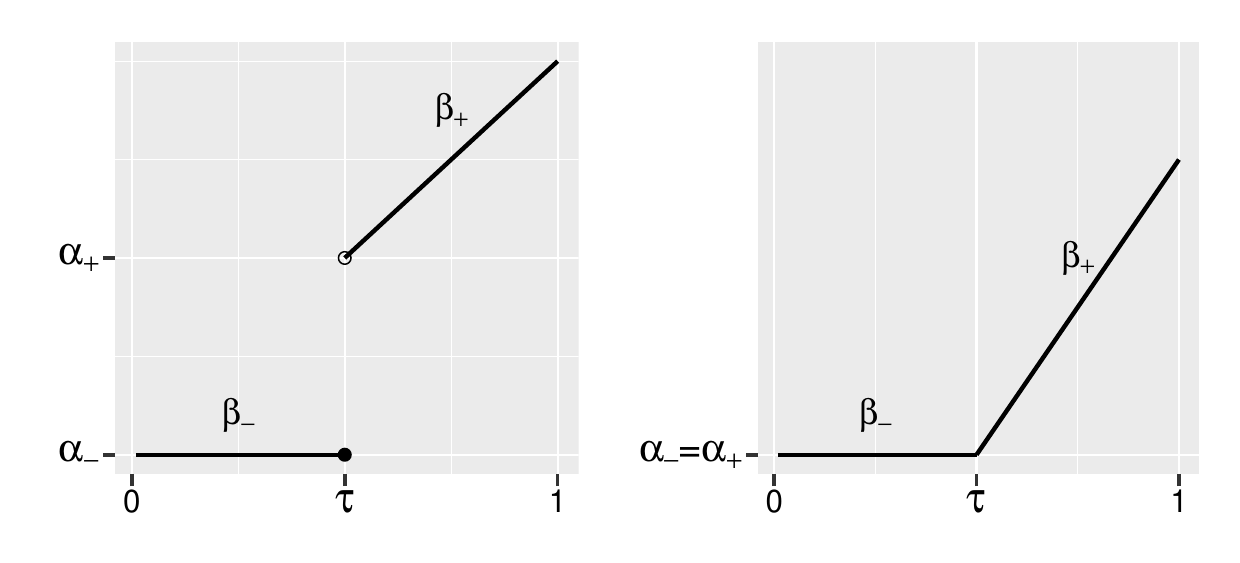}
    \caption{Sketches of some segmented linear functions \(f_\theta\) with a jump (left) and a kink (right).}
    \label{fig:model}
\end{figure}

In this paper, we consider the sequential detection of a change in the segmented linear regression model, where observations \(X_1,\ldots,X_n\) are given by 
\begin{subequations}
    \begin{equation}\label{eq:model}
X_i \coloneqq f_{\theta}\bigl(\frac{i}{n}\bigr)+\sigma \varepsilon_i,\qquad i=1,\ldots,n,  
\end{equation}
with \(\varepsilon_i \iidsim \mathcal{N}\l(0,1\r)\), the standard Gaussian distribution. For simplicity, we assume that the standard deviation \(\sigma\) is known and set to $1$, though it can be $\sqrt{n}$-consistently preestimated from data (e.g.\ via local differences; see e.g.\ \citealp{HaMa90} and \citealp{DeMuWa98}) without affecting our results. The unknown function \(f_\theta:\l[0,1\r]\mapsto \mathbb{R}\) is piecewise linear and takes the form of
\begin{equation}\label{eq:signal}
f_{\theta}\bigl(\frac{i}{n}\bigr) \coloneqq \begin{cases}\beta_-\l(\frac{i}{n}-\tau\r) + \alpha_-& \text{if }\frac{i}{n} \le \tau, \\\beta_+\l(\frac{i}{n}-\tau\r)+ \alpha_+ &  \text{if }\frac{i}{n} > \tau,
\end{cases}
\end{equation}
where \(\theta = \l(\tau,\alpha_-,\alpha_+,\beta_-,\beta_+\r)\) lies in the parameter space 
\begin{equation}\label{eq:parameterspace}
\Theta_{\delta_0} \coloneqq \bigl\{\theta = \l(\tau,\alpha_-,\alpha_+,\beta_-,\beta_+\r) : \delta_0 \le \tau \le 1-\delta_0,\max\left(|\alpha_+ - \alpha_-|,\,|\beta_+ -\beta_-|\right)\ge \delta_0\bigr\},
\end{equation}
\end{subequations}
for some (maybe unknown) \(\delta_0\in\l(0,{1}/{2}\r)\). A key aspect of our analysis is the distinction between two types of  structural changes in \(f_\theta\): a \emph{jump} and a \emph{kink}. We classify a change in \(f_\theta\) as a {jump} if \(\l|\alpha_+-\alpha_-\r|\ge \delta_0\), indicating a (significant) discontinuity in function values, and as a {kink} if \(\l|\alpha_+-\alpha_-\r|< \delta_0\) but \(\l|\beta_+-\beta_-\r|\ge \delta_0\), indicating a (significant) discontinuity in slope without a (significant) jump in values, see \Cref{fig:model}. The respective parameter spaces are denoted by \(\Theta^J_{\delta_0} \coloneqq \l\{\theta\in\Theta_{\delta_0} : |\alpha_+-\alpha_-|\ge \delta_0\r\}\) for the jump case and \(\Theta^K_{\delta_0} \coloneqq \Theta_{\delta_0}\setminus \Theta^J_{\delta_0}\) for the kink case.

Formally, we define an \emph{online detector} \(\hat\tau_n \) as a Markov stopping time with respect to the filtration \(\l\{\mathcal{F}_t,1\le t\le n\r\}\), where \(\mathcal{F}_t\) is the \(\sigma\)-algebra generated by observations \(X_1,\ldots,X_t\). Let  \(\mathcal{T}_n\) denote the collection of all such online detectors.  Following \citet[Chapter~6]{korostelev}, we aim at achieving the convergence rate of the minimax quadratic risk of detection, as the sample size $n$ goes to infinity. Namely, we study the rate at which the quadratic risk
\[
R_n^*(\delta_0) \;  \coloneqq\;\inf_{\hat\tau_n\in\mathcal{T}_n}\sup_{\theta = (\tau,\alpha_{-},\alpha_{+},\beta_{-},\beta_{+})\in\Theta_{\delta_0}}\E_{\tau}\left[\left(\hat\tau_n-\tau\right)^2\right]
\]
converges to zero. This risk naturally balances type~I (false detection) and type~II (missed detection) errors, since it can be decomposed as
\begin{equation}\label{eq:typeIandtypeIIerror}
\E_{\tau}\left[\left(\hat\tau_n-\tau\right)^2\right] = \E_{\tau}\left[\left(\hat\tau_n-\tau\right)^2\mathbb{I}\l(\hat\tau_n < \tau\r)\right] + \E_{\tau}\left[\left(\hat\tau_n-\tau\right)^2\mathbb{I}\l(\hat\tau_n \ge \tau\r)\right].
\end{equation}
This quadratic risk framework, rooted in in-fill asymptotics, differs slightly from conventional frameworks in sequential change point detection (see \Cref{ss:literature} below), aligning more closely with \emph{offline} change point analysis (see e.g.\ \citealp{FrMuSi14,NHZ16,TOV20,ChKi21}), where all observations are available upfront. By framing the online problem in a manner analogous to offline analysis, our approach facilitates a rigorous comparison between online and offline setups. {\rev Furthermore, the second term in \eqref{eq:typeIandtypeIIerror}, together with the false alarm probability, implies an upper bound on the expected detection delay as follows
\begin{equation}\label{eq:detdelbound}
   \E_{\tau}\left[\left(\hat\tau_n-\tau\right)^2\mathbb{I}\l(\hat\tau_n \ge \tau\r)\right]\ge \bigl(\E_\tau\l[\hat\tau_n-\tau\mid\hat\tau_n\ge\tau\r]\bigr)^2 \bigl(1 - \P_{\tau}(\hat\tau_n< \tau)\bigr), 
\end{equation}
thereby connecting the quadratic risk to standard performance metrics in sequential change point detection (see \Cref{corollary:detection delay}).} In addition, we put emphasis simultaneously on the statistical optimality and the computational and memory efficiency. Our contributions are threefold: 
\begin{enumerate}[label=\roman*.]
    \item 
    We introduce the Fast Limited-memory Optimal Change (FLOC) detector, which operates in constant time per observation and requires only constant memory, particularly independent of the sample size (see \Cref{p:comp}). The implementation of FLOC is provided in an R package available on GitHub under \url{https://github.com/AnnikaHueselitz/FLOC}.
    \item\label{i:ctr:b} 
    We establish that the proposed FLOC detector achieves minimax optimal rates in quadratic risk for estimating the change point in segmented linear regression (see \Cref{theorem:upperbound,theorem:lowerbound}). {\rev Moreover, it reliably distinguishes between jump and kink types of structural changes with high probability (see \Cref{theorem:falsedetectiontypebound}).} 
   \item
    Our theoretical analysis reveals a fundamental \emph{phase transition} between jump and kink scenarios (see~\Cref{tab:rates}). Specifically, the minimax rate scales as \(\log (n)/n\) for a jump and as \((\log (n)/n)^{1/3}\) for a kink, underscoring the different levels of detectability between abrupt and more gradual structural changes in the regression function. 
\end{enumerate}

\begin{table}[h]
    \centering
    \caption{Phase transition between a jump and a kink in segmented linear regression for online (this paper) and offline \citep{Che21} setups. \label{tab:rates}}{
    \begin{tabular}{lll|ll|ll} \toprule
    & \multicolumn{2}{l|}{Optimal rate in quadratic risk}  & \multicolumn{2}{l|}{Computational complexity} & \multicolumn{2}{l}{Memory complexity}\\
           & online & offline & online & offline & online & offline\\ \midrule
          Jump  & \(\mathcal{O}\l(\log (n)/n\r)\) & \(\mathcal{O}(1/n)\)  & \(\mathcal{O}(1)\) per observation & $\mathcal{O}(n^2)$  & \(\mathcal{O}(1) \) & \(\mathcal{O}(n)\)\\
         Kink  & \(\mathcal{O}\bigl((\log (n)/n)^{1/3}\bigr)\) & \(\mathcal{O}\l((1/n)^{1/3}\r)\)  & \(\mathcal{O}(1)\) per observation & $\mathcal{O}(n^2)$ & \(\mathcal{O}(1)\) &\(\mathcal{O}(n)\)\\ \bottomrule
    \end{tabular}}
\end{table}

This phase transition between a jump and a kink has been similarly observed in offline change point estimation in segmented linear models \citep{Che21}, see also \citet{GoTsZe06} and \citet{ShHaHa22} for general regression models. The literature on offline change point estimation in segmented linear models is notably extensive (see e.g.\ \citealp{BaPe98}, \citealp{muggeo2003estimating}, \citealp{FrHoMu14}, and more recently, \citealp{LeSeSh16}, \citealp{cho2024detection} and \citealp{cho2025covariance} for high-dimensional extensions). In \Cref{tab:rates}, we summarize the phase transition scenarios in segmented linear models, comparing online and offline setups in terms of optimal rates,  computational efficiency and memory requirements. It is worth emphasizing that, unlike offline setups, where model validity is typically required globally over the whole domain of the function, our online detector FLOC is effective even with only local model validity around the change point.  This flexibility makes our approach particularly well-suited for applications where model structures are reliable only within specific localized regions, as demonstrated with the real-world data set discussed in \Cref{section:realdata}.

\subsection{Related work}\label{ss:literature}
\paragraph{Minimax optimality}
We review the literature that is most pertinent to our approach from an asymptotic minimax perspective. For an infinite sequence of observations $X_1, X_2, \ldots$, the statistical performance of an online detector $\tilde\tau$, taking values in $\mathbb{N}$, is often assessed by 
$$
\sup_{1 \le \kappa < \infty} \mathrm{D}_{\kappa}(\tilde\tau) \qquad\text{subject to }\; \E \left[\tilde\tau\mid\text{no change}\right] \ge \gamma
$$
in the asymptotic regime of $\gamma \to \infty$. Here, $D_{\kappa}(\cdot)$ measures the detection delay when the change occurs at $\kappa$, and  specific examples include
\begin{subnumcases}{\mathrm{D}_{\kappa}(\tilde\tau) = }
\mathrm{ess\,sup}\, \E[\max(\tilde \tau - \kappa +1,0)\mid X_1,\ldots,X_{\kappa-1};\, \text{change at }\kappa], & \text{see \citet{Lord71}, or}\label{e:Lv}\\
\E[\tilde \tau - \kappa +1\mid \tilde\tau > \kappa;\, \text{change at }\kappa], & \text{ see \citet{Poll85}.}\label{e:Pv}
\end{subnumcases}
Within this framework, \citet{Yao93} investigated the optimal detection of a jump, with respect to \eqref{e:Lv}, in linear regression models, and \citet{YaKrPo99} studied the optimal detection of a kink with respect to \eqref{e:Pv}. One can relate such results to the model \eqref{eq:model}--\eqref{eq:parameterspace} in this paper by the correspondence of $\hat\tau = \tilde{\tau}/n$ and $n \approx \gamma$, see \Cref{corollary:detection delay} for further details. 

Minimax optimality results have also been established beyond linear regression models. For instance, asymptotically optimal detection methods have been developed for mean shifts \citep{AuHo04, YuMPWaRi23}, non-parametric frameworks \citep{Che19, HoKoWa21}, and high-dimensional settings \citep{Chan17, ChWaSa22}. 

\paragraph{Computational efficiency} Despite their statistical optimality, the aforementioned methods will be confronted with computational challenges in large-scale data sets, as both their runtime and memory usage scale linearly with the number of cumulatively observed data samples for each incoming sample. Recently, in the submodel of \eqref{eq:model}--\eqref{eq:parameterspace} with piecewise constant signals (i.e.\ $\beta_{-} = \beta_{+} = 0$), \citet{RoEcFeRi23} provided an efficient algorithm via functional pruning for computing the online detector by \citet{YuMPWaRi23}, achieving $\mathcal{O}(\log(t))$ computational complexity and an average memory requirement of $\mathcal{O}(\log(t))$ at the time of $t$ observations, for $t = 1, \ldots, n$. Variants of it using amortized cost or information from previous iterates seem to achieve $\mathcal{O}(1)$ computation complexity per observation, as indicated by empirical evidence \citep{WaRoEcFe24}. In addition, \citet{ChWaSa22} introduced an online method for a high-dimensional mean shift detection with a computational and memory complexity per observation that is independent of the number of previous samples (i.e.\ constant in observed sample sizes). 

As an addition to the existing literature, the proposed online detector FLOC achieves constant computational and memory complexity per observation, while also attaining minimax optimal rates in terms of quadratic risk for segmented linear regression models simultaneously for a jump and a kink (see again \Cref{tab:rates}). 

\subsection{Organization and notation}

The remainder of the paper is organized as follows. In \Cref{section:upperbound}, we introduce formaly the FLOC detector and establish its statistical guarantees in estimating the change point. The minimax lower bounds on convergence rates in quadratic risk are given in \Cref{s:lb}. \Cref{section:numexp} examines the empirical performance of the proposed FLOC detector on both simulated and real-world data sets and in comparison to state-of-the-art methods. All the proofs are deferred to \Cref{appendix:proofs}, and \Cref{s:dis} concludes the paper with a discussion. 

Throughout, we denote by \(\P_{\tau}\l(\cdot\r)\) the probability measure and  by \(\E_{\tau}\l[\cdot\r]\) the expectation under the model \eqref{eq:model}--\eqref{eq:parameterspace} with the true change point~\(\tau\). In particular, \(\P_{0}\l(\cdot\r)\) and \(\E_{0}\l[\cdot\r]\) refer to the case with no change (i.e.\ $\tau = 0$). For sequences $\{a_n\}$ and $\{b_n\}$ of positive numbers, we use $a_n = \mathcal{O}(b_n)$ to indicate that $a_n \le C b_n$ for some finite constant $C > 0$, and write $a_n \asymp b_n$ if $a_n = \mathcal{O}(b_n)$ and $b_n = \mathcal{O}(a_n)$.

\section{FLOC: Statistical theory and computational complexity}\label{section:upperbound}

In this section, we introduce the Fast Limited-memory Optimal Change (FLOC) detector, based on two online detectors, one of which is used to detect a jump and the other to detect a kink. We then derive upper bounds on the quadratic risk of FLOC as well as its computational cost.  

\subsection{Risk bounds}\label{section:riskbounds}

\begin{figure}[ht]
    \centering
    \includegraphics[width = .8\linewidth]{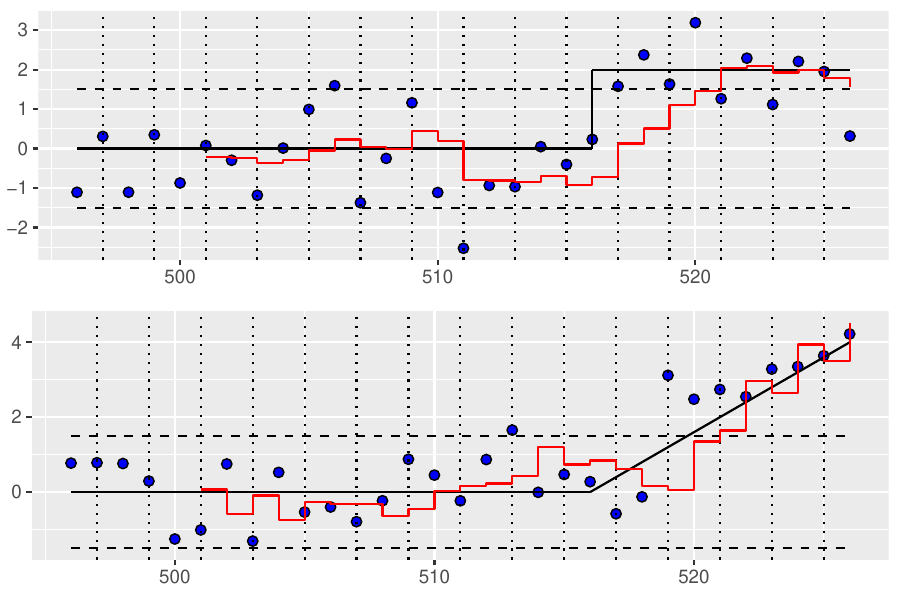}
    \caption{Simulated scenarios illustrating a jump (upper panel) and a kink (lower panel). The black lines represent the underlying signal \(f_\theta\), which exhibits a change at  \(\tau = 516/n\), and the blue dots denote observations sampled from a Gaussian distribution with mean \(f_\theta\), as specified in \eqref{eq:model}.   The red step functions show the values of the test statistic computed by FLOC. The vertical dotted lines indicate observation points where FLOC discards a bin; the bin sizes are \(N_J = N_K = 2\). The pre-change signal is estimated using the first \(k = 500\) observations. A change is flagged when the test statistic exceeds the threshold, with \(\rho_J = 1.5\) and  \(\rho_K = 0.3\), shown as horizontal dashed lines.  For both the jump and the kink scenarios, this occurs at observation $521$.  For visual clarity in the lower panel (kink case), the threshold and test statistic values are scaled by a factor of 5.}
    \label{fig:algorithm_sample}
\end{figure}

The core idea of the online detection algorithm is to estimate the pre-change signal \(f_-\) at the beginning of the data sequence and then continuously assess whether the data within the most recent time window remains consistent with this underlying function.  The online detector is built upon a certain test statistic, which reports a change if the test statistic exceeds a pre-specified threshold.  To detect a jump, the test statistic is the mean of the residuals between the observed data and the estimated pre-change signal over the latest interval. For a kink, a weighted mean of these residuals is used instead, assigning greater weights to more recent observations.  

\begin{subequations}

Consider a data sequence \(X_1,X_2,\ldots, X_n\), and let $k, N_J, N_K \le n$ be natural numbers that may depend on \(n\). A minimax rate optimal choice of these parameters will be provided in \Cref{theorem:upperbound}, while practical selection strategies are discussed later in \Cref{section:parameterchoice}. The function $f_-$ is estimated using the first $k$ data points via least squares, as follows:
\begin{equation}\label{def:hatfminus}
    \hat f_-\bigl(\frac{i}{n}\bigr) \coloneqq \hat \alpha + \hat \beta \frac{i}{n},\qquad\text{ with }\quad \bigl(\hat \alpha, \hat \beta\bigr) \coloneqq \underset{\l(\alpha, \beta\r)\in\mathbb{R}^2}{\arg\min} \sum_{i=1}^{k} \left(\alpha + \beta \frac{i}{n} - X_i\right)^2.
\end{equation}
We define the \emph{CUSUM jump test statistic} of bin size $N_J$ as
\begin{equation}\label{def:jumpstat}
    J_m \coloneqq \frac{1}{M_J}\sum_{i=1}^{M_J}\l(X_{m-M_J+i}-\hat f_-\bigl(\frac{m-M_J+i}{n}\bigr)\r), \qquad k < m\le n,
\end{equation}
with \(M_J \coloneqq 2N_J +\l( m \mod N_J\r)\), and similarly the \emph{CUSUM kink test statistic} of bin size $N_K$ as
\begin{equation}\label{def:kinkstat}
    K_m \coloneqq \frac{6}{M_K\l(M_K+1\r)\l(2M_K+1\r)}\sum_{i=1}^{M_K}i\l(X_{m-M_K+i}-\hat f_-\bigl(\frac{m-M_K+i}{n}\bigr)\r),\qquad k<m\le n,
\end{equation}
with \(M_K \coloneqq 2N_K +\l( m \mod N_K\r)\). {\rev For an intuitive understanding of \(J_m\) and \(K_m\), we observe that both test statistics are proportional to the partial derivatives with respect to $\alpha$ and $\beta$ of the residual sum of squares as in \eqref{def:hatfminus} evaluated at \((\hat\alpha,\hat\beta)\) , respectively, making them sensitive to changes in these parameters. Specifically, when there is no change in $\alpha$ (or $\beta$), the corresponding test statistic $J_m$ (or $K_m$) remains close to zero, whereas it increases in magnitude when a change occurs.} Thus, we introduce, for some thresholds \(\rho_J, \rho_K > 0\), the corresponding \emph{online detectors} $\hat \tau_{J,n}$ and $\hat \tau_{K,n}$ as the first time the test statistics exceed their respective thresholds, i.e., 
\begin{equation}\label{def:jumpdetector}
\hat \tau_{J,n} \coloneqq \begin{cases}1 &\text{if } |J_m| < \rho_J\text{ for all  }m\in\{k+1,\ldots,n\},\\
\min\bigl\{m: |J_m| \ge \rho_J,\ k< m \le n\bigr\}/n&\text{otherwise,}
\end{cases}
\end{equation}
for a jump and 
\begin{equation}\label{def:kinkdetector}
    \hat \tau_{K,n} \coloneqq \begin{cases}1 &\text{if } |K_m| < \rho_K\text{ for all  }m\in\{k+1,\ldots,n\},\\
    \min\bigl\{m: |K_m| \ge \rho_K,\ k < m \le n\bigr\}/n&\text{otherwise,}
    \end{cases}
\end{equation}
for a kink.
Then we define the \emph{FLOC (Fast Limited-memory Optimal Change) detector} as the minimum of the above two detectors, i.e.,
\begin{equation}\label{def:FLOC}
    \hat\tau_n \coloneqq \min\l(\hat \tau_{J,n}, \hat \tau_{K,n}\r).
\end{equation}
\end{subequations}
It is straightforward to verify that the detectors $\hat \tau_{J,n}$ and $\hat \tau_{K,n}$, as defined in \eqref{def:jumpdetector} and \eqref{def:kinkdetector}, are Markov stopping times. As a consequence, FLOC, being the minimum of two Markov stopping times, is also a Markov stopping time, that is, FLOC is indeed an online detector (see also the proof in \Cref{appendix:proofs}). {\rev An illustration of FLOC is provided in \Cref{fig:algorithm_sample}. While the bin sizes $N_J$ and $N_K$ are fixed, the window sizes used by both the CUSUM jump and kink test statistics in \eqref{def:jumpstat} and \eqref{def:kinkstat} vary over time. This design offers an advantage in terms of storage efficiency: instead of retaining all recent observations and removing them individually, it suffices to store only the weighted sums of observations and discard entire bins when starting a new one (see \Cref{ss:com}).}


\begin{theorem}\label[theorem]{theorem:upperbound}
    Assume the segmented linear regression model \eqref{eq:model}--\eqref{eq:parameterspace}, and let $\hat f_-$ be defined in \eqref{def:hatfminus} with $k = c n$ for some constant $c \in (0,\delta_0)$. Then:
    \begin{enumerate}[label=\roman*.]
        \item Let \(\hat\tau_{J,n}\) be the CUSUM jump detector in \eqref{def:jumpstat} and \eqref{def:jumpdetector} with bin size $N_J = 10^3 \log (n)/(2 c^2) $ and threshold $\rho_J = 4c/5$. Then, it holds for sufficiently large $n$,
        \[
        \sup_{\theta = (\tau,\alpha_{-},\alpha_{+},\beta_{-},\beta_{+})\in\Theta^J_{\delta_0}}\mathbb{E}_{\tau}\left[\left(\frac{n}{\log n}\l(\hat\tau_{J,n}-\tau\r)\right)^2\right]\;\le\; r^*_J < \infty,
        \]
        where \(r^*_J\) is independent of $n$ and can be chosen as \(r^*_J = 9\cdot 10^6\cdot 4^{-1}\delta_0^{-4} + 1 \).
        \item Let \(\hat\tau_{K,n}\) be the CUSUM kink detector in \eqref{def:kinkstat} and \eqref{def:kinkdetector} with bin size $N_K = (300/c^2)^{1/3} n^{2/3}\log^{1/3}(n)$ and threshold $\rho_K = 4c/(5n)$. Then, it holds for sufficiently large $n$,
        \[
        \sup_{\theta = (\tau,\alpha_{-},\alpha_{+},\beta_{-},\beta_{+})\in\Theta^K_{\delta_0}}\mathbb{E}_{\tau}\left[\left(\frac{n^{1/3}}{\log^{1/3}n}\l(\hat\tau_{K,n}-\tau\r)\right)^2\right]\le r^*_K < \infty,
        \]
        where \(r^*_K\) is independent of $n$ and can be chosen as \(r^*_K =9 \cdot 300^{2/3}\delta_0^{-4/3} + 1\).
        \item For the FLOC detector  $\hat\tau_n$ defined in \eqref{def:FLOC}, both of the above rates hold for sufficiently large $n$.
    \end{enumerate}    
\end{theorem}

The idea of the proof (for details see \Cref{appendix:proofs}) is to characterize an event on which the detector reports a change point without false alarms and achieves a detection delay that is in the range of the considered rate; then such an event is shown to occur with a high probability. The constants in \Cref{theorem:upperbound} are explicit, but they are not optimized and could potentially be improved.

{\rev 
The originally assumed requirement of $k \asymp n$ historical data points can be substantially relaxed. As shown in the proof in \Cref{appendix:proofs}, it suffices to choose $k \asymp n^{2/3}\log(n)^{1/3}$ in the jump case, and $k \asymp n^{8/9}\log(n)^{1/9}$ in the kink case. In the special case of a piecewise constant model, i.e., when $\beta_-=\beta_+ = 0$ in \eqref{eq:signal}, the requirement can be further reduced to $k \asymp \log(n)$. The gap between the requirements $k \asymp n^{2/3}\log(n)^{1/3}$ and $k \asymp \log(n)$ reveals the additional sample complexity of historical data incurred in the more general piecewise linear setting. 

This requirement on $k$ has practical implications when applying the online detector to monitor multiple change points. In online settings, the method can be naturally extended to the multiple change point scenario by restarting the detector each time a change is declared. In such a framework, the minimum segment length between two change points must be at least of order $k$. 
}

As demonstrated later in \Cref{s:lb}, the rates established in \Cref{theorem:upperbound} are matched by corresponding lower bounds, revealing that $\hat\tau_{J,n}$ is minimax optimal for detecting a jump, while $\hat\tau_{K,n}$ is minimax optimal for detecting a kink. Further, FLOC, defined as the minimum of these two, is minimax optimal in detection of the change point, simultaneously for jump and kink types of structural changes. {\rev Motivated by this observation, one might naturally consider inferring the type of structural change, either a jump or a kink, by examining which of the two detectors $\hat\tau_{J,n}$ and $\hat\tau_{K,n}$ triggers an alarm first. The following proposition establishes that this approach is indeed theoretically justified.

\begin{proposition}\label[proposition]{theorem:falsedetectiontypebound}
    Assume the same setup as in \Cref{theorem:upperbound} and let $n$ be sufficiently large. Then:
        \begin{align*}
        & \P_\tau\l(\hat\tau_{K,n}<\hat\tau_{J,n}\r) \le  2n^{-3}\qquad\text{if}\quad\theta=(\tau,\alpha_{-},\alpha_{+},\beta_{-},\beta_{+}) \in \Theta^J_{\delta_0},\\
       \text{and}\quad & \P_\tau\l(\hat\tau_{J,n}<\hat\tau_{K,n}\r)\le  2n^{-3}\qquad\text{if}\quad\theta= (\tau,\alpha_{-},\alpha_{+},\beta_{-},\beta_{+}) \in \Theta^K_{\delta_0}.
       \end{align*}
\end{proposition}

The proof (see \Cref{appendix:proofs}) builds on the high probability events introduced in the proof of \Cref{theorem:upperbound}, under which the false alarm rates of both detectors are well controlled. In the presence of a jump, the kink detector reacts more slowly due to its larger bin size. Conversely, when a kink is present, the induced change in slope is too gradual to trigger the jump detector before the kink detector detects the change. The rate of upper bound $n^{-3}$ could be replaced by $n^{-C}$ for any constant $C \ge 0$, while doing so would require enlarging the bin sizes, which would in turn worsen the constants $r_J^*$ and $r^*_K$ in \Cref{theorem:upperbound}. This result confirms that each detector is very unlikely to falsely signal first in the incorrect regime. 

Furthermore, using the decomposition in \eqref{eq:typeIandtypeIIerror}, Theorem~\ref{theorem:upperbound} yields statistical guarantees on the false alarm probability, the average run length and the expected detection delay (cf.~\Cref{ss:literature}).

\begin{subequations}
\begin{corollary}\label[corollary]{corollary:detection delay}
        Assume the same setup as in \Cref{theorem:upperbound} and let $n$ be sufficiently large. Then:
    \begin{enumerate}[label=\roman*.]
        \item For the jump detector \(\hat\tau_{J,n}\), it holds that \(\P_\tau(\hat\tau_{J,n} < \tau) \le n^{-3}\), 
        \(\E_0\l[\hat\tau_{J,n}\r]\ge 1 - n^{-3}\)
        and
        \begin{equation}\label{e:delay:jump}
        \sup_{\theta = (\tau,\alpha_{-},\alpha_{+},\beta_{-},\beta_{+})\in\Theta^J_{\delta_0}} \frac{n}{\log n} \E_\tau\l[\hat\tau_{J,n}-\tau\mid\hat\tau_{J,n}\ge\tau\r]\;\le\; \sqrt{r^*_J} < \infty.
        \end{equation}
        \item For the kink detector \(\hat\tau_{K,n}\), it holds that  \(\P_\tau(\hat\tau_{K,n} < \tau) \le n^{-4}\), 
        \(\E_0\l[\hat\tau_{K,n}\r]\ge 1 - n^{-4}\)  and
        \begin{equation}\label{e:delay:kink}
        \sup_{\theta = (\tau,\alpha_{-},\alpha_{+},\beta_{-},\beta_{+})\in\Theta^K_{\delta_0}}\frac{n^{1/3}}{\log^{1/3} n} \E_\tau\l[\hat\tau_{K,n}-\tau\mid\hat\tau_{K,n}\ge\tau\r]\;\le\; \sqrt{r^*_K} < \infty.
        \end{equation}
        \item For the FLOC detector $\hat\tau_n$, we have   \(\P_\tau(\hat\tau_{n} < \tau) \le n^{-3}+n^{-4}\) and \(\E_0\l[\hat\tau_{n}\r]\ge 1 - n^{-3} - n^{-4}\), and the same upper bounds as in \eqref{e:delay:jump} and \eqref{e:delay:kink} hold.  
    \end{enumerate}
\end{corollary}
\end{subequations}

As noted for \Cref{theorem:falsedetectiontypebound}, the bounds of $n^{-3}$ and $n^{-4}$ may be replaced by
$n^{-C}$ for any constant $C\ge0$, provided appropriate parameter values are chosen. In this regard, the parameters can be tuned to achieve a user-specified false alarm probability or average run length, see \Cref{section:parameterchoice} for details. }

{\rev
In addition, it is worth noting that, in \Cref{theorem:upperbound}, \Cref{theorem:falsedetectiontypebound} and \Cref{corollary:detection delay}, the Gaussian assumption on the noise is used solely to control the deviation of weighted sums of the noise variables. More precisely, in the proofs in \Cref{appendix:proofs}, we employ the concentration inequality
\[
\mathbb{P} \left(\left\vert\sum_{i = 1}^n w_i \varepsilon_i \right\vert \ge t \right) \le \exp\left(-\frac{t^2}{2}\right)\qquad \text{for}\quad t > 0,
\]
where the weights $w_i$ are deterministic and satisfy $\sum_{i = 1}^n w_i^2 = 1$. Similar concentration bounds hold for sums of independent sub-Weibull random variables, see e.g.\ \citet{KuCh22}. Further, certain forms of temporal dependence can be incorporated using the functional dependence framework (introduced in \citealp{Wu05}). As a consequence, the established statistical guarantees of our proposed online detectors can be extended to these more general noise settings.
}

\subsection{Computational complexity}\label{ss:com}
A pseudocode implementation of FLOC is provided in \Cref{algorithm:floc}. The algorithm maintains sums and weighted sums in bins of a fixed size to calculate the test statistics. When the third bin is full, the first bin is discarded, and a new bin is initiated. The benefits of varying window sizes over a fixed window size become evident in the implementation of FLOC. While a constant window size achieves the same asymptotic properties, it would require storing individual observations to allow for gradual updates as the window moves. As the window size depends on \(n\), the storage cost for a fixed window size would grow with the rate of observations. In contrast,  \Cref{algorithm:floc} achieves computational and memory costs of $\mathcal{O}(1)$, independent of the window size and the number of observation until the current time. This property is particularly advantageous for online detectors, as minimal computational and memory requirements are often considered essential design criteria (cf.\ Introduction).
\begin{proposition}\label[proposition]{p:comp}
    For the FLOC detector as defined in \eqref{def:jumpstat}--\eqref{def:FLOC}, the following holds:
    \begin{enumerate}[label=\roman*.]
        \item The computational cost is \(\mathcal{O}(1)\) per incoming observation.
        \item The storage cost is \(\mathcal{O}(1)\).
    \end{enumerate}
\end{proposition}
{\rev We note that a full batch computation of $\hat f_{-}$ in \eqref{def:hatfminus} using the historical data of size $k$ requires $\mathcal{O}(k)$ runtime, but this is performed only once. Alternatively, $\hat f_{-}$ can be updated sequentially, resulting in an $\mathcal{O}(1)$ computational cost per data point.}  

\begin{algorithm}
\caption{FLOC: Fast Limited-memory Optimal Change detector}\label{algorithm:floc}
\begin{algorithmic}[1]
    \Require{jump sums \(S_{J,1}, S_{J,2}, S_{J,3}\), kink sums \(S_{K,1}, S_{K,2}, S_{K,3}\), weighted sums \(W_1, W_2, W_3\)}
    \Ensure{latest data point \(X_{t}\) at time \(t\), jump bin size \(N_J\), kink bin size \(N_K\), jump threshold \(\rho_J\), kink threshold \(\rho_K\), estimate of pre-change signal \(\hat f_-\) }
    \State compute \(r_J\) \(\gets t \mod N_J\)
    \State compute \(r_K\) \(\gets t \mod N_K\)
    \State compute scaling factor \(d\) \(\gets {(2N_K+r_K+1)(2N_K+r_K+2)(4N_K+2r_K+3)}/{6}\)
    \If{\(r_J = 0\)}
    \State update jump sums \(S_{J,1} \gets S_{J,2}\), \(S_{J,2}\gets S_{J,3}\), \(S_{J,3} \gets 0\)
    \EndIf
    \If{\(r_K = 0\)}
    \State update kink sums \(S_{K,1} \gets S_{K,2}\), \(S_{K,2}\gets S_{K,3}\), \(S_{K,3} \gets 0\)
    \State update weighted sums \(W_1 \gets W_2\), \(W_2\gets W_3\), \(W_3 \gets 0\)
    \EndIf
    \State update jump sum  \(S_{J,3} \gets S_{J,3} + X_t - \hat f_-\l(t\r)\)
    \State update kink sum  \(S_{K,3} \gets S_{K,3} + X_t - \hat f_-\l(t\r)\)
    \State update weighted sum \(W_3 \gets W_3 + \l(r+1\r)\bigl(X_t - \hat f_-\l(t\r)\bigr)\)
    \State compute the CUSUM jump test statistic \(J\gets \l(S_{J,1}+S_{J,2}+S_{J,3}\r)/\l(2N_J+r_J+1\r)\)
    \State compute the CUSUM  kink test statistic \(K\gets  \l(W_1 + W_2 + W_3 + N_K S_{K,2} + 2N_K S_{K,3}\r) / d\)
    \If{\(|J|\ge \rho_J\)}
    \State \Return{jump detected}
    \ElsIf{\(|K|\ge \rho_K\)}
    \State \Return{kink detected}
    \Else
    \State \Return{no change detected}
    \EndIf
\end{algorithmic}
\end{algorithm}

\section{Lower risk bounds}\label{s:lb}
In this section, we complement the upper bounds in the previous section with matching lower bounds, and discuss the relation to existing results in the literature. 

The following theorem provides a lower bound for all online detectors (i.e.\ Markov stopping times). 
\begin{theorem}\label[theorem]{theorem:lowerbound}
        Assume the segmented linear regression model \eqref{eq:model}--\eqref{eq:parameterspace}. Then, there exist positive constants \(r_{*J}\) and \(r_{*K}\), independent of \(n\), such that
    \begin{align*}
    &\liminf_{n\to\infty}\inf_{\hat{\tau}_n\in \mathcal{T}}\max_{\theta = (\tau,\alpha_{-},\alpha_{+},\beta_{-},\beta_{+})\in\Theta^J_{\delta_0}} \mathbb{E}_{\tau}\left[\left(\frac{n}{\log n}\l(\hat{\tau}_n-\tau\r)\right)^2\right] \ge r_{*J} > 0,\\
    \text{and}\quad& 
 \liminf_{n\to\infty}\inf_{\hat{\tau}_n\in \mathcal{T}}\max_{\theta = (\tau,\alpha_{-},\alpha_{+},\beta_{-},\beta_{+})\in\Theta^K_{\delta_0}} \mathbb{E}_{\tau}\left[\left(\frac{n^{1/3}}{\log^{1/3}n}\l(\hat{\tau}_n-\tau\r)\right)^2\right] \ge r_{*K} >0.\end{align*}
\end{theorem}
The proof of \Cref{theorem:lowerbound}, detailed in \Cref{appendix:proofs}, proceeds by contradiction.  In the proof, we utilize a simplification by fixing the changes in jump and slope, as the lower bound will remain valid for the full parameter space. Then we consider potential change points within the unit interval, spaced in proportion to the target rate. If the specified rate were not a valid lower bound, there would be a detector such that the probability of correct detection tends to $1$. However, by analyzing the likelihood ratio between two distinct change points and leveraging the property of Markov stopping times, we show that the number of possible change points grows too fast for this to occur.


In combination with \Cref{theorem:upperbound}, \Cref{theorem:lowerbound} shows that the two distinct rates for online detection of jumps and kinks are minimax optimal. These rates closely align with the offline setup, differing only in log factors (see \Cref{tab:rates} in the Introduction). Notably, FLOC attains the minimax optimality in detection of both jumps and kinks. For the detection of a jump, the optimal rate for the segmented linear regression model coincides with the optimal rate for the submodel with piecewise constant signal and a single jump (see e.g.\ Chapter~6 in \citealp{korostelev}). It reveals that asymptotically a change in the slope is negligible in comparison to a jump. Further, for the submodel with piecewise constant signal, \citet{YuMPWaRi23} showed a rate of order \(\log\frac{n}{\alpha}\) on the detection delay under the condition \(\P\l(\hat\tau< \infty \mid \text{no change}\r)\le \alpha\), which ensures that under the null hypothesis the false alarm probability of online detector is at most  \(\alpha\). The computational cost for their proposed detector scales linearly with the number of  observations available at time $t$, and this complexity was later improved to $\mathcal{O}(\log t)$ using functional pruning \citep{RoEcFeRi23}. In contrast, our FLOC detector, designed for a more general model,  attains the same detection delay rate in this specific submodel, and requires a computational and memory costs independent of time $t$, see also \Cref{ss:literature} and \Cref{corollary:detection delay}.

\section{Numerical experiments}\label{section:numexp}

In this section, we first introduce an empirical method for the selection of  thresholds in FLOC (\Cref{algorithm:floc}). Subsequently, we assess the empirical performance of FLOC on simulated datasets by benchmarking it against state-of-the-art online change point detection methods. In addition, we investigate the robustness of FLOC in settings that go beyond the scope of our theoretical guarantees, particularly in the presence of non-Gaussian noise. Finally, we demonstrate the application of FLOC to the analysis of COVID-19 excess mortality data from the United States. The implementation used in this section is publicly available as an R package on GitHub at \url{https://github.com/AnnikaHueselitz/FLOC}.

{\rev
\subsection{Choice of parameters}
\label{section:parameterchoice}

FLOC involves five tuning parameters. \Cref{theorem:upperbound} provides theoretical guidance for their selection: the historical data size  \(k=cn\), with \(c < \delta_0\) (i.e., smaller than the minimal detectable change), thresholds \(\rho_J = 4c/5\) and \(\rho_K= 4c/(5n)\) and bin sizes \(N_J = 10^3\log(n)/(2c^2)\) and \(N_K = (300/c^2)^{1/3} n^{2/3} \log^{1/3}(n)\). However, these settings may be impractical in cases where  \(\delta_0\) is unknown. We therefore provide empirical recommendations below, based on simulation studies.

The amount of historical data $k$ governs the accuracy of estimating the pre-change distribution. Larger values of $k$ yield more reliable estimates and result in a more stable detector. See \Cref{tab:performance} for simulation results illustrating the effect of $k$ on detection delay. In practice, one should utilize as much historical data as is available or computationally feasible. When $k$ is small, higher thresholds are needed to control the type I error, which in turn increases the detection delay. We can also incorporate a known pre-change signal $f_-$ by directly substituting it in place of its estimate $\hat f_-$.
}

In general, increasing the detection threshold reduces the typeI error at the cost of increasing the typeII error. The type~I error is typically measured by the \emph{false alarm probability} \(\P_{\tau}\l(\hat\tau < \tau\r)\), which is the probability that a false alarm, when the change point is at \(\tau\), or by the \emph{average run length} under the null hypothesis \(\E_0\l[\hat\tau\r]\), which is the expected run length if no change occurs. 
The type II error is often measured by the \emph{expected detection delay} \(\E_{\tau}\l[\hat\tau-\tau\mid\hat\tau \ge \tau\r]\). 
In applications, \(n\) does not need to be specified beforehand. Thus, we will use in the following the expected detection delay and the average run length in terms of the number of observations.

\begin{figure}
    \centering
    \includegraphics[width = .85\linewidth]{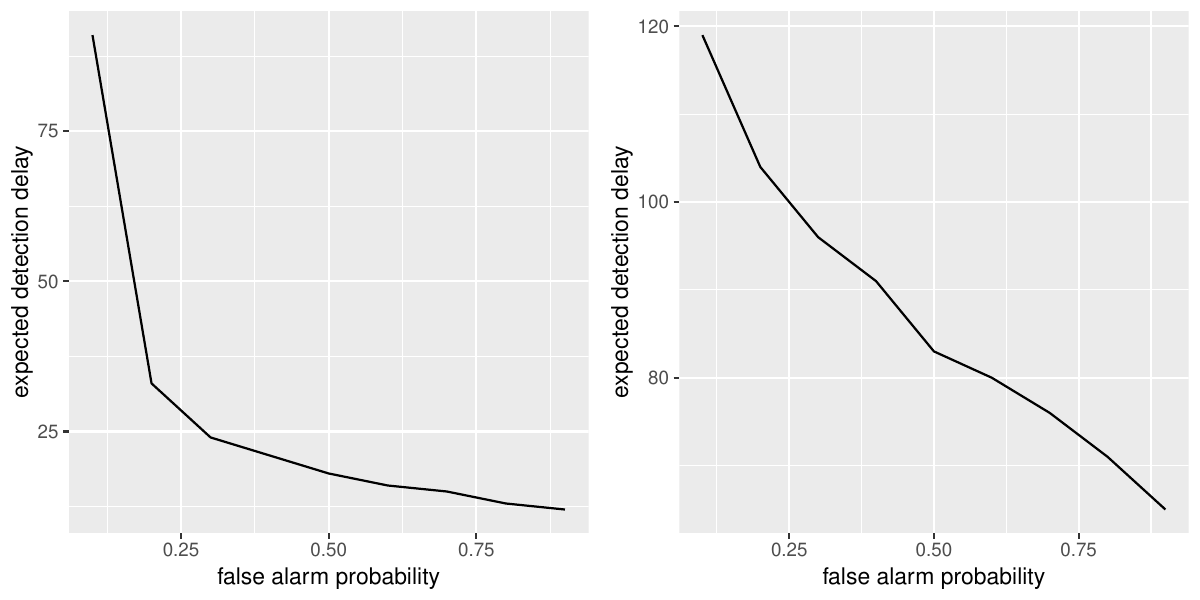}
    \caption{Comparison of the type I error, expressed as the false alarm probability \(\P(n\hat\tau_n < 1000)\), and the type II error, measured by the expected detection delay for a change point at observation \(500\), for the jump detection (left) and the kink detection (right) using FLOC. In both plots, the bin size is \(N_J = N_K = 5\) and \(f_-\) is estimated using \(k = 500\) observations. Threshold tuning is performed with \(r = 5000\) data sets and the expected detection delay is estimated with \(200\) repetitions. The jump size used is \(1\) and the change in slope is \(0.01\). }
    \label{fig:typeIvstypeIIerror}
\end{figure}

\begin{figure}
    \centering
    \includegraphics[width = .85\linewidth]{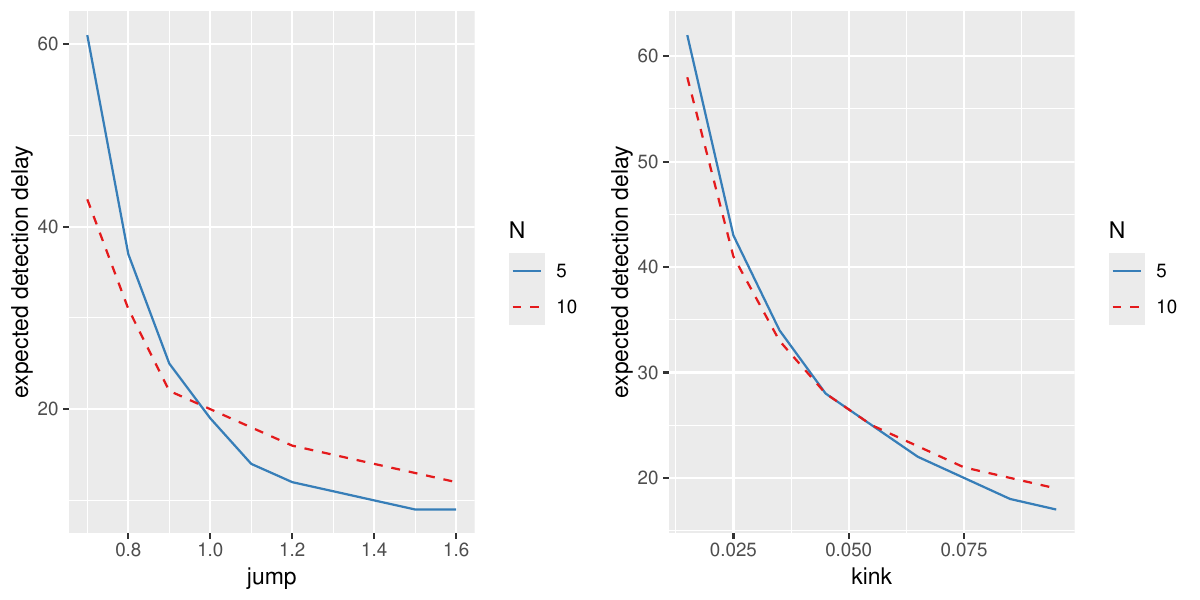}
    \caption{Expected detection delay when the true change point occurs at observation 
$500$, for bin sizes of 5 (blue solid) and 10 (red dashed), for a jump (left) and a kink (right). The thresholds are calibrated to achieve a false alarm probability of  \(\P\l(n\hat\tau_n\le 1000\r)\approx 0.5\). The pre-change signal $f_-$ is estimated using $k = 500$ observations. For threshold calibration, $r = 1000$ data sets are generated, and the expected detection delay is estimated using $500$ independent repetitions.}
    \label{fig:binsizecomparison}
\end{figure}

\begin{figure}
    \centering
    \includegraphics[width = .85\linewidth]{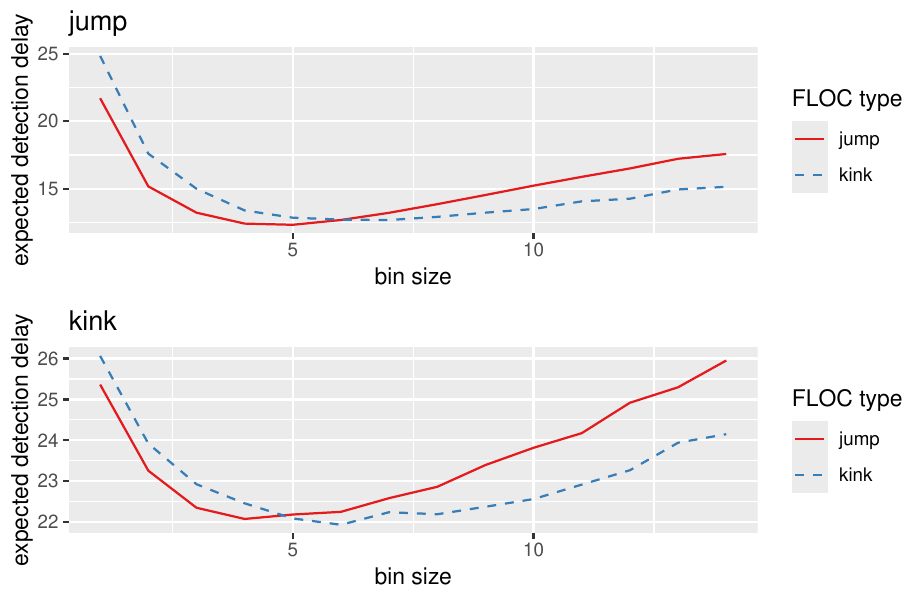}
    \caption{\rev Comparison of the detection of a jump (red solid) and a kink (blue dashed) using FLOC across multiple bin sizes \(N_J\) and \(N_K\) for a jump of size \(1\) (top) and a kink of size \(0.05\) (bottom). Thresholds are tuned to achieve an average run length of \(500\), and the size of historical dat is \(k = 500\). The expected detection delay is estimated over \(1000\) repetitions.}
    \label{fig:floctype}
\end{figure}

\Cref{fig:typeIvstypeIIerror} illustrates the relationship between the false alarm probability and the expected detection delay.  
In practice, thresholds can be tuned to control the type~I error at a user-specified level \(\eta\) e.g.~by bounding the false alarm probability. 
This amounts to finding the smallest threshold such that \(\P_{\tau}\l(\hat\tau<\tau\r)\le \eta\) for some fixed~\(\tau\). More precisely, for given bin sizes \(N_J\) and \(N_K\) and \(k\) historical observations, we generate independently $r$  data sets \(X_1^{(m)},\ldots,X_{k+n\tau}^{(m)}\) for \(1 \le m\le r\) under the null hypothesis. For each data set, we compute the maximum statistics \(X^*_{J,m} \coloneqq \max_{1\le i <n\tau} J_i^{(m)} \) for jump detection and \(X^*_{K,m} \coloneqq \max_{1\le i <n\tau} K_i^{(m)} \) for kink detection. 
{\rev
 By default, FLOC simultaneously monitors for both a jump and a kink; however, it can be configured to detect only one of these changes by setting the threshold $\rho_K$ (or $\rho_J$) to infinity, thereby disabling the detection of kinks (or jumps).}
When using FLOC to detect only jumps 
, we select the jump threshold $\rho_J$ as the \(\left(1-\eta\right)\)-quantile of \(\{X^*_{J,m}: 1\le m \le r\}\). 
A similar procedure applies when detecting only kinks. When FLOC is applied to monitor both change types simultaneously, the thresholds can be selected as appropriate quantiles of \(\{X^*_{J,m}: 1\le m \le r\}\)  and \(\{X^*_{K,m}: 1\le m \le r\}\) such that \(\P_{\tau}\l(\hat \tau_J < \tau\r) \approx \P_{\tau}\l(\hat \tau_k < \tau\r)\) and \(\P_{\tau}\l(\hat \tau_J < \tau, \hat \tau_k < \tau\r) \approx \eta\). The effectiveness of this approach is examined by simulation in \Cref{tab:performance}.
{\rev
A related strategy can be employed to calibrate thresholds with respect to the \emph{average run length}. Due to the rolling window structure, FLOC exhibits approximate memorylessness under the null, suggesting that the run length distribution is roughly exponential. To target an average run length of at least  \(\E_0\l[\hat\tau\r]\ge \tau\), one can use the same threshold selection procedure as for the false alarm probability, with  \(\eta = 1- 1/e\). \Cref{tab:arl} reports thresholds obtained via this approach along with the  empirical average run lengths that are actually achieved.

\begin{table}
    \centering
    \caption{Performance of FLOC (\Cref{algorithm:floc}) on simulated data. The thresholds $\rho_J$ and $\rho_K$ are calibrated with \(10000\) repetitions to achieve a target false alarm probability (FA) of $0.5$, when the change point is at observation \(1000\) for $k\in \{500,1000\}$ historical observations, and the change point at \(10000\) for \(k = 5000\). The actual false alarm probability and the expected detection delay in terms of number of observations (last six columns), under a change point at \(k\)-th observation, are estimated via Monte Carlo simulations with 200 repetitions, across varying magnitudes of jumps and kinks.}
    \begin{tabular}{l c c c c c c c c c c c c}\toprule
        & & & & & & &\multicolumn{3}{c}{jump}&\multicolumn{3}{c}{kink}\\ \cmidrule(rl){8-10}\cmidrule(rl){11-13}
        FLOC & \(N\) & target FA & $k$ & \(\rho_J\) & \(\rho_K\) & FA & 2 & 1 & 0.5 & 0.5 & 0.1 & 0.02\\ \midrule
jump  & 5     & 0.5  & 500   & 1.031 & -     & 0.5   & 7     & 21    & 329   & -     & -     & -    \\ 
  & 10    & 0.5  & 500   & 0.749 & -     & 0.46  & 10    & 19    & 159   & -     & -     & -    \\ 
  & 10    & 0.5  & 1000  & 0.658 & -     & 0.53  & 9     & 17    & 63    & -     & -     & -    \\ 
  & 10    & 0.5 & 5000  & 0.79  & -     & 0.47  & 11    & 20    & 127   & -     & -     & -    \\ 
  & 15    & 0.5  & 500   & 0.627 & -     & 0.55  & 13    & 24    & 171   & -     & -     & -    \\ \midrule
kink  & 5     & 0.5  & 500   & -     & 0.148 & 0.48  & -     & -     & -     & 6     & 17    & 52   \\ 
  & 10    & 0.5  & 500   & -     & 0.058 & 0.51  & -     & -     & -     & 8     & 18    & 47   \\ 
  & 10    & 0.5  & 1000  & -     & 0.051 & 0.54  & -     & -     & -     & 7     & 17    & 44   \\ 
  & 10    & 0.5 & 5000  & -     & 0.062 & 0.52  & -     & -     & -     & 9     & 19    & 52   \\ 
  & 15    & 0.5  & 500   & -     & 0.033 & 0.46  & -     & -     & -     & 8     & 18    & 49   \\ \midrule
both  & 5     & 0.5  & 500   & 1.071 & 0.15  & 0.5   & 6     & 17    & 352   & 6     & 16    & 49   \\ 
  & 10    & 0.5  & 500   & 0.781 & 0.06  & 0.54  & 9     & 19    & 221   & 8     & 19    & 47   \\ 
 & 10    & 0.5  & 1000  & 0.687 & 0.05  & 0.46  & 8     & 15    & 73    & 8     & 17    & 43   \\ 
  & 10    & 0.5 & 5000  & 0.818 & 0.06  & 0.46  & 10    & 20    & 133   & 9     & 20    & 52   \\ 
  & 15    & 0.5  & 500   & 0.651 & 0.03  & 0.54  & 11    & 22    & 192   & 8     & 18    & 48   \\ 
 \bottomrule
    \end{tabular}
    \label{tab:performance}
\end{table}

\begin{table}
    \caption{{\rev Thresholds \(\rho_K\) and \(\rho_J\), tuned using the procedure outlined in \Cref{section:parameterchoice} to achieve target average run lengths (ARLs), along with the corresponding actual ARLs. Results are reported for the detection of a jump, a kink, and both of them, across different bin sizes  \(N \in \{10, 15\}\) and historical data sizes \(k \in\{ 1000, 2500, 5000\}\). Threshold calibration is performed using 10000 repetitions, and ARLs and expected detection delays (last six columns) are estimated from 100 independent repetitions.}}
    \centering
    \begin{tabular}{l c c c c c c c c c c c c c}
        \toprule
        & & & & & & &\multicolumn{3}{c}{jump}&\multicolumn{3}{c}{kink}\\ \cmidrule(rl){8-10}\cmidrule(rl){11-13} 
        FLOC & \(N\) & target ARL & \(k\) & \(\rho_J\) & \(\rho_K\) & ARL & 2& 1& 0.5 & 0.5& 0.1 &0.02\\
        \midrule
        jump & 10 & 1000 & 1000 & 0.621 & - &  987.73&  9 & 16 &  54 & - & - & -\\ 
        & 15 & 1000 & 1000 & 0.497 & - & 1064.14& 10 & 19 &  41 & - & - & -\\
        & 10 & 1000 & 5000 & 0.584 & - & 1087.70&  8 & 15 &  33 &- & - & -\\
        &10 & 5000 & 2500 & 0.734 & - & 4617.44 & 10 & 19 &  96 & - & - & -\\  
        \midrule
        kink & 10 & 1000 & 1000 & - & 0.0487 &  932.22&  - & - &  - & 7 & 15 & 42\\
        & 15 & 1000 & 1000 & - & 0.0267 &  999.05&  - & - &  - & 7 & 17 & 41\\
        & 10 & 1000 & 5000 & - & 0.046 &  922.33&  - & - &  - & 7 & 16 & 42\\ 
        & 10 & 5000 & 2500 & - & 0.057 & 4682.73 & - & - & - & 8 & 19 & 49\\ 
        \midrule
        both & 10 & 1000 & 1000 & 0.65 & 0.0509 &  975.99&  7 & 13 &  48 & 7 & 17 & 40\\  
        & 15 & 1000 & 1000 & 0.522 & 0.0278 &  962.65&  8 & 16 &  43 & 7 & 17 & 42\\
        & 10 & 1000 & 5000 & 0.612 & 0.048 &  941.56&  7 & 13 &  39 & 7 & 16 & 39\\
        & 10 & 5000 & 2500 & 0.763 & 0.06 & 5244.13 & 9 & 17 & 112 & 8 & 19 & 47\\
        \bottomrule
    \end{tabular}
    \label{tab:arl}
\end{table}

For the bin sizes \(N_J\) and \(N_K\), smaller values are preferable for detecting abrupt large changes, whereas larger values are more suitable for identifying gradual small changes, see \Cref{fig:binsizecomparison} for an illustration. Consequently, an optimal bin size exists for each given change magnitude. This is demonstrated in \Cref{fig:floctype}, which contrasts the expected detection delay across a range of bin sizes. The figure also shows that, with the optimal bin size (i.e.~the one that minimizes the expected detection delay), the jump detector in FLOC is the first to signal a change when the underlying change is a jump, and analogously for the kink detector. This observation aligns with the theoretical guarantee in \Cref{theorem:falsedetectiontypebound}. Moreover, the flatness of the performance curves around the optimal bin size indicates that FLOC is relatively robust to the choice of bin sizes. In practice, we therefore recommend fine-tuning the bin sizes based on the anticipated type and magnitude of the change in a given application. If such prior knowledge is unavailable, we suggest applying FLOC with a pair of bin sizes, specifically, a small and a large bin size. This combination enables detection of both abrupt and gradual changes (see \Cref{fig:mult_bins}) and, compared with a single intermediate bin size, delivers more reliable performance across a broad range of change magnitudes. Including additional bin sizes offers no significant improvement but increases computational cost.}

\begin{figure}
    \centering
    \includegraphics[width = .85\linewidth]{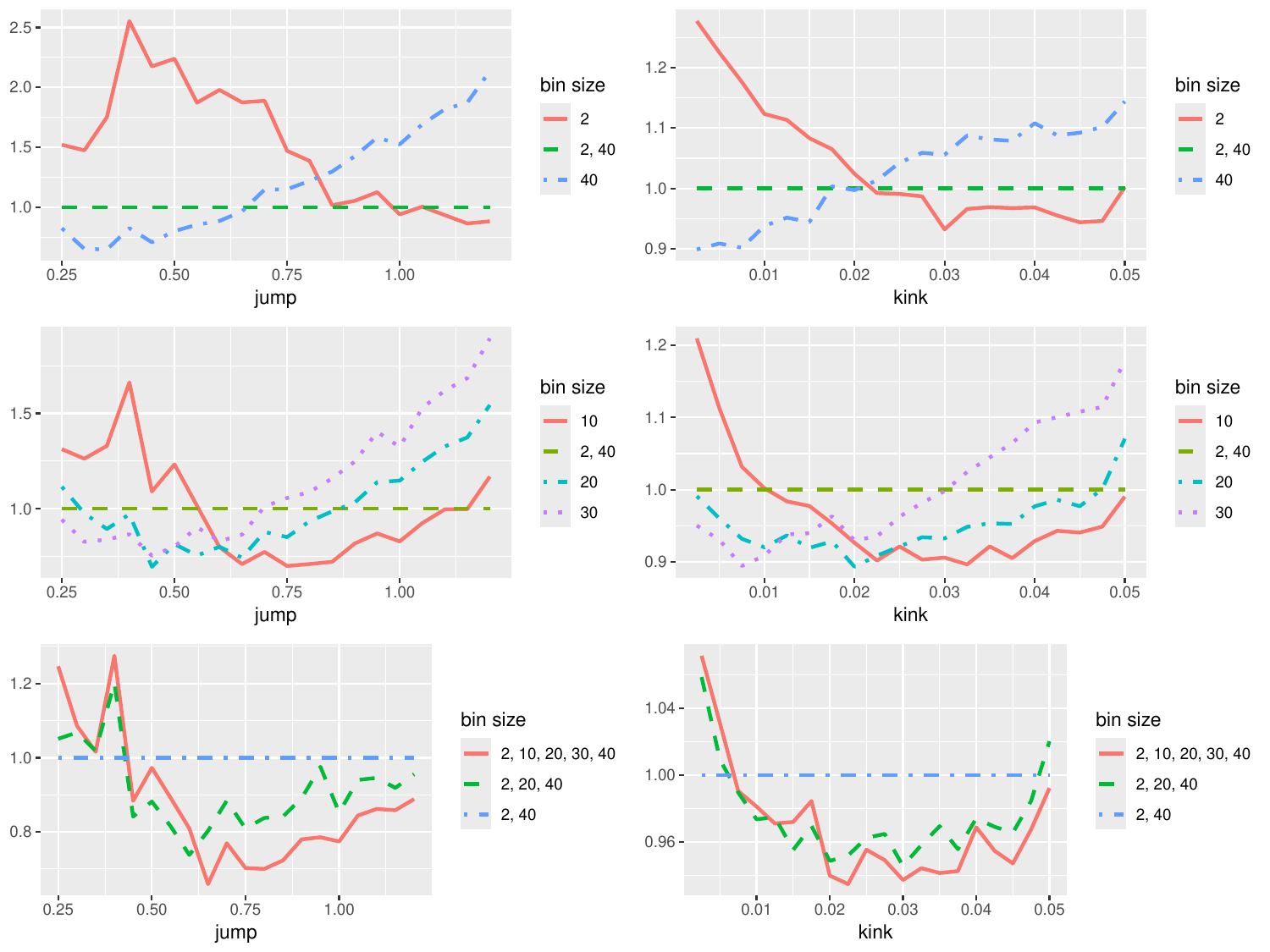}
    \caption{\rev Relative performance of FLOC with different bin size choices, measured by the expected detection delay across various change magnitudes for a jump (left) and a kink (right). The baseline is FLOC with bin sizes \(\{2, 40\}\). The vertical axis reports the ratio of expected detection delay relative to the baseline. For all bin size choices, thresholds are calibrated to achieve an average run length of \(500\) observations with \(k = 500\) historical data points. Expected detection delays are computed over \(500\) repetitions.}
    \label{fig:mult_bins}
\end{figure}

{\rev
\subsection{Comparison study}

We benchmark the performance of FLOC against two state-of-the-art methods, which we denote by FOCuS (Functional Online CuSum; \citealp{RoEcFeRi23}) and Yu-CUSUM \citep{YuMPWaRi23}. We use the  implementations of both methods available on GitHub (\url{https://github.com/gtromano/FOCuS} and \url{https://github.com/HaotianXu/changepoints}). FLOC is designed to detect structural changes involving both jumps and kinks under a segmented linear regression model, whereas FOCuS and Yu-CUSUM are tailored for detecting mean shifts in a piecewise constant mean model. To ensure a fair comparison, we restrict our evaluation to scenarios with piecewise constant means, i.e., we set  \(\beta_{-} = \beta_{+} = 0\) in \eqref{eq:signal}. All three methods are provided with the same amount of historical data, with $k \in \{1000, 2000, 2500, 5000\}$. The data are standardized using the standard deviation of the historical segment, and each method is calibrated to achieve a target average run length of either \(1000\) or \(5000\).

We report expected detection delays under varying jump sizes in \Cref{tab:comparison}. The results show that FOCuS generally performs among the best across most scenarios, while FLOC demonstrates comparable or slightly superior performance when the bin size is appropriately chosen. In contrast, Yu-CUSUM consistently exhibits inferior performance.

As actual runtimes can vary due to differences in implementation, we provide a theoretical comparison of the computational complexity per observation. Yu-CUSUM requires either $\mathcal{O}(n^2)$ runtime and  $\mathcal{O}(1)$ storage, or $\mathcal{O}(n)$ runtime and $\mathcal{O}(n)$ storage. FOCuS achieves an expected runtime and storage complexity of $\mathcal{O}(\log n)$. In contrast, FLOC offers the desirable constant runtime and storage complexity, i.e.~$\mathcal{O}(1)$.

\begin{table}
    \caption{
    {\rev Expected detection delays of FLOC, FOCuS \citep{RoEcFeRi23} and Yu-CUSUM \citep{YuMPWaRi23} under a piecewise constant mean model of varying jump sizes. All methods are tuned to achieve a common average run length (ARL) and employ  \(k\) historical observations from the pre-change distribution. The notation FLOC \(N\) indicates the choice of bin size \(N \in \{5, 10, 15, 20\}\). The change point is fixed at 100 observations. Expected detection delays are estimated over 100 independent repetitions, with the lowest value in each setting highlighted in bold.} }
    \centering
    \begin{tabular}{l c c c c c c c c}
        \toprule
        ARL & \(k\) & jump & FLOC 5 & FLOC 10 & FLOC 15 & FLOC 20 & FOCuS & Yu-CUSUM\\ 
        \midrule
        1000 & 1000 & 1  & 12.58 &  14.39 &  17.96 &  20.32 & \textbf{10.14} & 13.67 \\ 
         &  & 0.5  &  104.51 &  52.29 &  60.13 &  49.25&  \textbf{43.53} & 58.88\\ 
         & 2000 & 1  &   12.34 &  14.75 &  17.32 &  22.03&  \textbf{11.95} & 14.93\\ 
         &  & 0.5 &   61.61 &  43.04 &  \textbf{39.76} &  \textbf{39.76}&  \textbf{39.76} & 51.04\\ 
         & 5000& 0.5  &   56.48 &  36.21 &  38.58 &  \textbf{35.36}&  36.96 &    49.38\\ 
         &  & 0.25  &  244.22 & 157.10 & \textbf{131.29} & 140.53& 148.27 &   238.26\\
        5000 & 2500 & 1  &   18.68 &  18.82 &  22.11 &  25.11&  \textbf{15.22} & 23.25\\ 
         &  & 0.5 &  137.92 &  82.80 &  72.47 &  66.65 &  \textbf{54.97} & 89.18\\ 
        \bottomrule
    \end{tabular}
    \label{tab:comparison}
\end{table}

\subsection{Non-Gaussian noise}

As a robustness study against model misspecification, we evaluate the performance of FLOC under non-Gaussian noise, specifically Student's \(t\)-distributions with varying degrees of freedom to model different tail behavior. All data are standardized using the empirical standard deviation of the historical observations. We consider two scenarios: detecting a jump in a piecewise constant model and detecting a kink in a segmented linear model. In both cases, FLOC is tuned under the (misspecified) assumption of Gaussian noise to achieve a target average run length of \(1000\) observations. The simulation results are summarized in \Cref{tab:t-dist}, where we include FOCuS as a competitor, also tuned to the same target average run length. 

For jump detection, FLOC demonstrates better robustness to heavy-tailed noise than FOCuS, consistently achieving higher average run length and shorter expected detection delays under all degrees of freedom. In particular, the performance of FLOC closely matches the Gaussian case when the degrees of freedom exceed \(4\).

For kink detection, FLOC shows comparable robustness: its performance is nearly indistinguishable from the Gaussian case when the degrees of freedom exceed \(2\). A distinction is that the actual average run length tends to be slightly below the target value in the kink scenario, whereas it exceeds the target in the jump scenario. This is likely due to differences in the weighting schemes of the respective test statistics. In comparison to FOCuS, FLOC has higher average run length, and shorter or similar expected detection delays. }

\begin{table}[h]
    \centering
        \caption{{\rev Comparison of FLOC and FOCuS \citep{RoEcFeRi23} under Student's \(t\)-distributed noise with varying degrees of freedom (DF). Both methods are calibrated to achieve a target average run length (ARL) of \(1000\) under the (misspecified) assumption of Gaussian noise. For FLOC, the bin sizes are set to  \(N_J = N_K = 15\). We use \(k = 5000\) historical observations from the pre-change distribution and standardize subsequent observations using the empirical standard deviation of the historical data. Both the ARL and the expected detection delay (EDD) are estimated based on \(500\) independent repetitions. We consider two scenarios of a jump of magnitude \(0.5\) and of a kink with a change in slope of \(0.01.\) The last row, with \(\text{DF} = \infty\), corresponds to the Gaussian case.}}
    \begin{tabular}{l r c r c r c r c}
        \toprule
        & \multicolumn{4}{c}{Jump} &  \multicolumn{4}{c}{Kink}\\ \cmidrule(lr){2-5}\cmidrule(lr){6-9}
       & \multicolumn{2}{c}{FLOC} &\multicolumn{2}{c}{FOCuS} & \multicolumn{2}{c}{FLOC} &\multicolumn{2}{c}{FOCuS}  \\ \cmidrule(lr){2-3}\cmidrule(lr){4-5} \cmidrule(lr){6-7}\cmidrule(lr){8-9}
        DF & ARL & EDD & ARL & EDD  & ARL & EDD & ARL & EDD\\ 
        \midrule
        1 &  3764 & 43  & 1707 & 51 &3035 & 68 & 1707 & 70\\ 
        2 &   2012 & 40 &176 & 42 & 1305 & 62 & 176 & 57\\ 
        3 &   1258 & 38 &128 & 39 &  647 & 59 & 128 & 50\\ 
        4 &    903 & 35  &141 & 39 &  735 & 59 & 141 & 54 \\ 
        5 &   1027 & 37 &184 & 42 &  866 & 60 & 184 & 55\\ 
        10 &  1068 & 39 &319 & 43 &   984 & 59 & 319 & 61\\ 
        30 &  1031 & 37 &698 & 41 &  1041 & 59 & 698 & 61\\ 
        \(\infty\) & 1089 & 37 & 960 & 44&   959 & 60  & 960 & 62\\
        \bottomrule
    \end{tabular}
    \label{tab:t-dist}
\end{table}

\subsection{Real data}\label{section:realdata}

\begin{figure}[!ht]
    \centering
    \includegraphics[width = 0.9\linewidth]{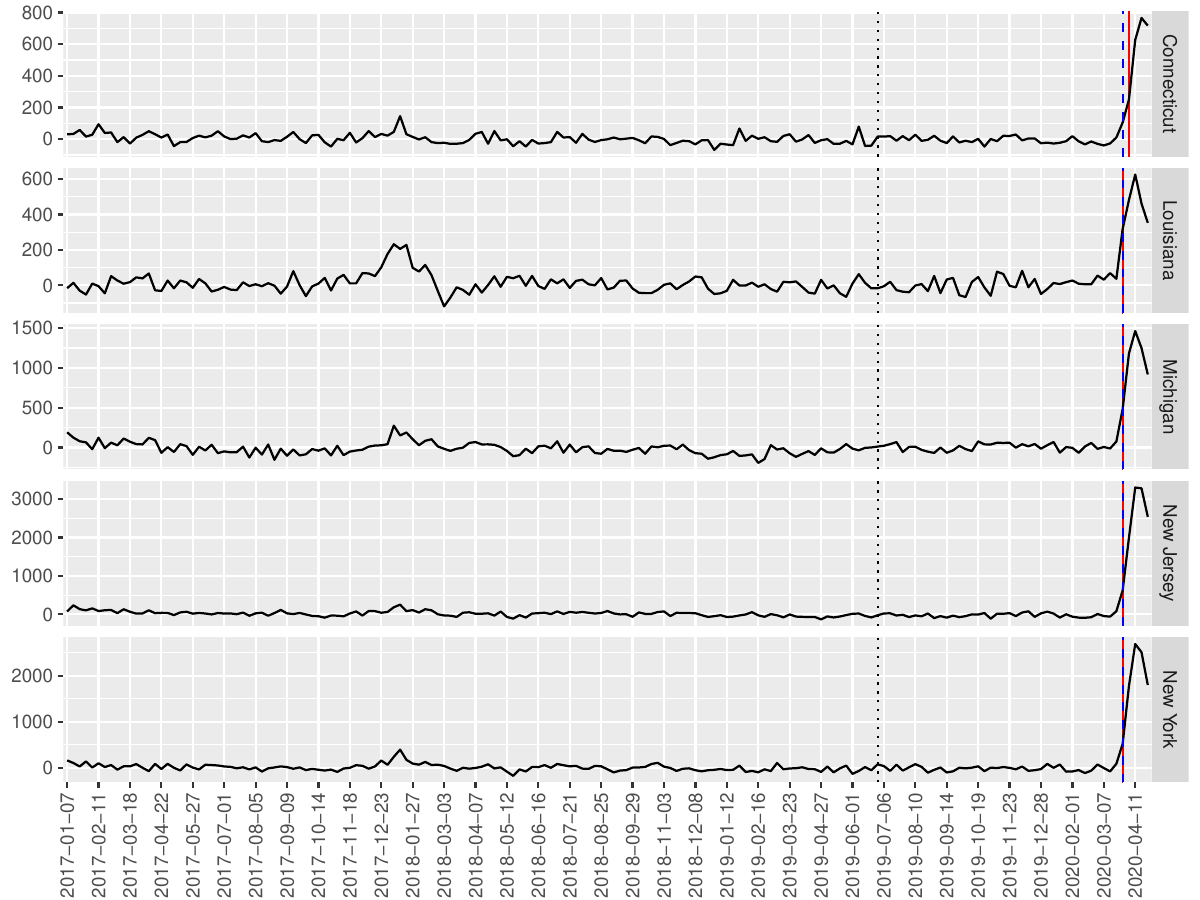}
    \caption{
    Excess death data in the states identified as significant by \texttt{ocd\_CI} \citep{chen2024inference}. The blue dashed line marks the change detected by \texttt{ocd\_CI}, while red lines indicate changes detected by FLOC. The data up to the dotted line (2020-06-30) serves as historical data for both approaches.}
    \label{fig:excessdeath}
\end{figure}

\begin{figure}[!ht]
    \centering
    \includegraphics[width = 0.9\linewidth]{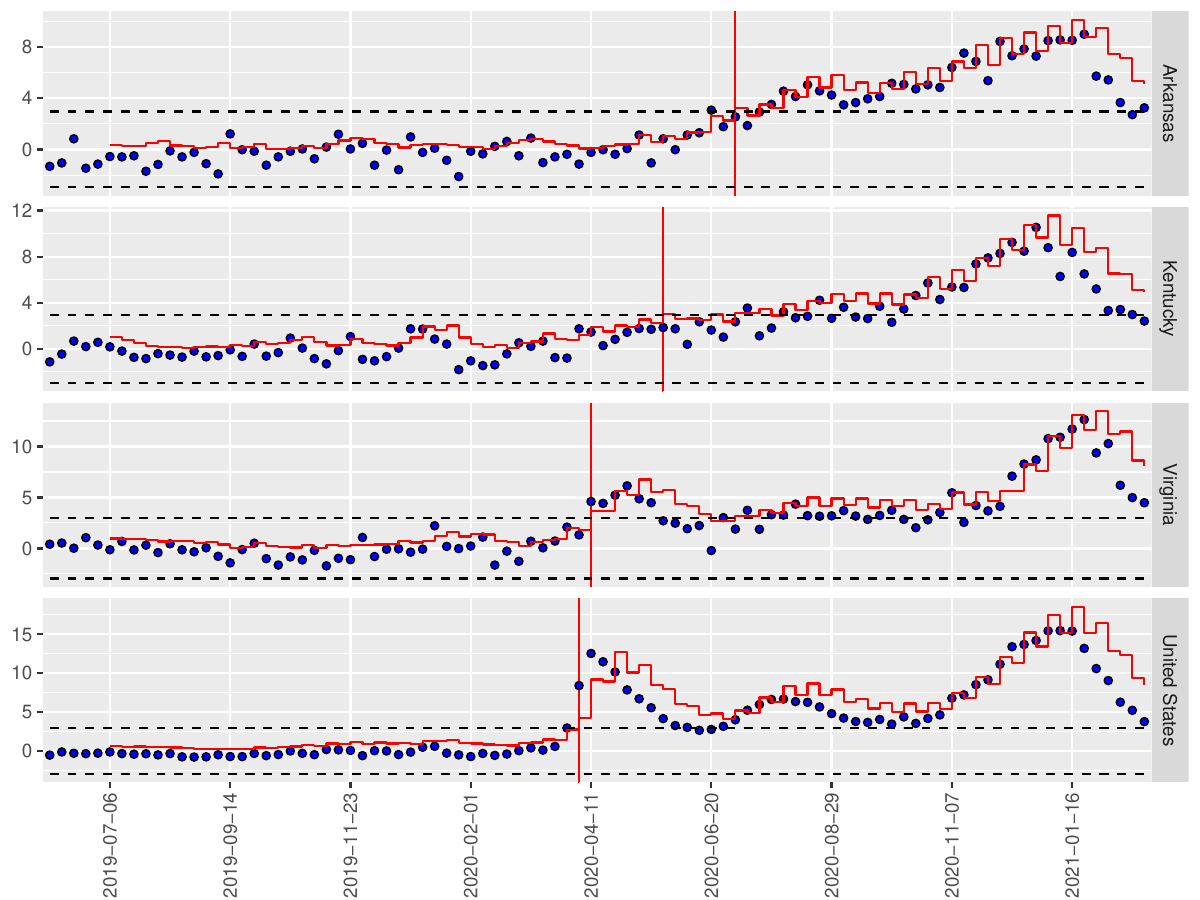}
    \caption{
    Standardized excess death data for Arkansas, Kentucky, Virginia, and the United States (pooled data). The red step function represents  the test statistic values of FLOC, scaled by a factor of \(4\) for visibility. The horizontal dashed line marks the detection threshold, and the vertical red line indicates the detected change in slope.}
    \label{fig:excessdeath_visual}
\end{figure}

To evaluate the empirical performance of FLOC on real-world data, we analyze excess mortality data from the CDC website \url{https://www.cdc.gov/nchs/nvss/vsrr/covid19/excess_deaths.htm}. This data set consists of estimates of excess deaths in the United States from the beginning of 2017 through September of 2023 by states (see \Cref{fig:excessdeath,fig:excessdeath_visual}), which were calculated based on weekly counts of deaths using Farrington surveillance algorithms (\citealp{NoEnFaGaAnCha13}, see the CDC website for further details). For our analysis, we treat the data through the end of June 2019 as historical data, and standardize the entire data set with the mean and the standard deviation of the historical period.  The primary objective of this analysis is to study the impact of Covid-19. This data set was also used by \citet{chen2024inference} to demonstrate the effectiveness of their online algorithm \texttt{ocd\_CI}, which detects changes in high-dimensional means. We thus include \texttt{ocd\_CI} as a baseline for comparison with FLOC.

The increase in excess deaths to the onset of Covid-19 appears to exhibit a gradual change in slope, rather than abrupt changes in values (see again \Cref{fig:excessdeath,fig:excessdeath_visual}). Therefore, we apply only the kink component of the FLOC algorithm by setting the jump threshold to infinity. We choose a bin size of \(N_K = 2\) weeks,  and tune the (kink) threshold \(\rho_K = 0.738\) to maintain a false alarm probability of approximately \(0.01\) on the remaining data.  The FLOC algorithm is applied individually to each state as well as to the pooled data for the United States. In contrast, the \texttt{ocd\_CI} method is applied to the entire data set, identifying a common change across all states and reporting a subset of states where this change occurs. 

\Cref{fig:excessdeath} shows the detection results of FLOC and \texttt{ocd\_CI} across the five states identified by \texttt{ocd\_CI} as having significant changes. In four of these five states (New York, New Jersey, Louisiana and Michigan), both algorithms detect changes during the week ending March 28, 2020. For the fifth state (Connecticut), FLOC identifies a change in the week ending April 4, 2020, which slightly differs from \texttt{ocd\_CI} (March 28, 2020). These detected changes align well with the timeline of Covid-19 outbreak in the United States as documented by the CDC (at \url{https://www.cdc.gov/museum/timeline/covid19.html}). In addition to these states, FLOC identifies changes in regions with a more moderate rise in excess deaths, where \texttt{ocd\_CI} does not detect any changes. For example, FLOC detects changes in Virginia during the week ending April 11, 2020, in Arkansas during the week ending July 4, 2020, and in Kentucky during the week ending May 23, 2020.  Further, FLOC detects a change in the pooled data for the United States in the week ending March 28, 2020. \Cref{fig:excessdeath_visual} presents a selection of the excess death curves and the corresponding changes detected by FLOC.

\begin{figure}[!tp]
    \centering
    \includegraphics[width =0.9 \linewidth]{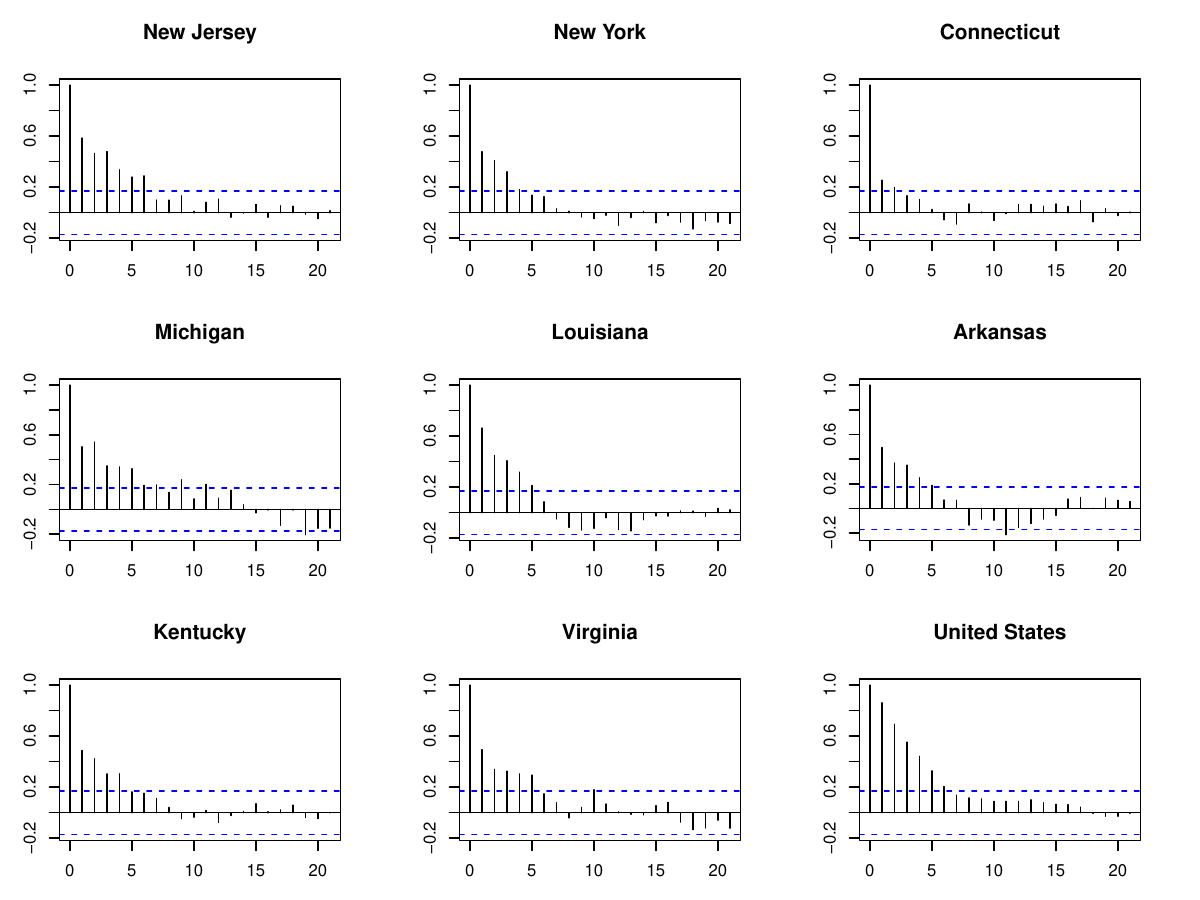}
    \caption{ACF plots showing the autocorrelation of excess death data (up to June 30, 2020) at various time lags for selected states and the pooled data for the United States.}
    \label{fig:excessdeath_acf}
\end{figure}

While FLOC yields insightful results on this data set, we acknowledge potential violations of the model assumptions  in \Cref{section:model}. For instance, the data exhibits positive correlation, as shown by the autocorrelation function (ACF) plots in \Cref{fig:excessdeath_acf}. We note that \texttt{ocd\_CD} also does not take the serial dependency into account (see \citealp{chen2024inference}). Furthermore, the underlying signal adheres to a kink-like structure only locally around the change points. A comprehensive analysis that incorporates the dependence structure of the data is beyond the scope of this paper and warrants further investigation.

\section{Proofs}\label{appendix:proofs}
\allowdisplaybreaks

\begin{proof}[Proof of \Cref{theorem:upperbound}]
The proof is divided in two parts. In the first part we show the upper bound for a jump at the change point and in the second part the upper bound for a kink at the change point. The statement for the FLOC algorithm then follows naturally, {\rev because 
\begin{align*}
\bigl(\min\{\hat\tau_{J,n}, \hat\tau_{K,n} \} - \tau\bigr)^2  & = \bigl(\min\{\hat\tau_{J,n}, \hat\tau_{K,n} \} - \tau\bigr)_+^2 + \bigl(\tau - \min\{\hat\tau_{J,n}, \hat\tau_{K,n} \} \bigr)_+^2\\
& =  \min_{\theta \in\{J,K\}}  \bigl(\hat\tau_{\theta,n} - \tau\bigr)_+^2 + \max_{\theta \in\{J,K\}}  \bigl(\tau - \hat\tau_{\theta,n}\bigr)_+^2 \\
& \le \max\left\{ \bigl(\hat\tau_{J,n} - \tau\bigr)^2, \, \bigl(\hat\tau_{K,n} - \tau\bigr)^2 \right\}
\end{align*}
where $(x)_+:=\max\{x, 0\}$.}
We denote by \(f_- (\cdot) \coloneqq \beta_-\l(\cdot-\tau\r) + \alpha_- \) and \(f_+(\cdot) \coloneqq \beta_+\l(\cdot-\tau\r) + \alpha_+\) the two linear functions that correspond to the two segments of \(f_{\theta_0}(\cdot)\).

\paragraph{For a jump} We assume that a jump occurs at the change point, namely, \(|\alpha_+-\alpha_-|\ge \delta_0\), for some \(\delta_0>0\). Our detector \(\hat\tau_{J,n}\) is given in \eqref{def:jumpdetector} with threshold \(\rho_J = 4c/5\) for some \(c \in (0,1)\). Recall that \(N_J\coloneqq b\log n\), with some \(b\) independent of \(n\), is the number of observations in a bin and  \(M_J = 2N_J +\l( m \mod N_J\r)\) is the window size at time \(m\), on which we calculate the test statistics. It follows that \(2N_J\le M_J < 3N_J\), which we assume to be integers for notational ease. We write only \(N_J = N\) and \(M_J = M\) for the jump part.

The detector \(\hat\tau_{J,n}\) is a stopping time, because \(\{\hat \tau_{J,n}= s/n\}  =  \{J_s \ge \rho_J,J_m < \rho_J,k\le m \le s-1\}\). Note that \(\{\hat\tau_{J,n} = s/n\}\) is determined by \(J_1,\ldots, J_s\), which are determined by \(X_1,\ldots,X_s\) and therefore independent of \(X_{s+1},\ldots,X_{n}\). 
Thus, \(\{\hat\tau_{J,n}= s/n\}\) is \(\mathcal{F}_s\)-measurable.

Denote by \(\tau\) the true change point.
Then, for \(m < n \tau\), i.e.\ before the true change point
\begin{subequations}
    \begin{align}
    J_m &= \left|\frac{1}{M}\sum_{i=1}^{M}\varepsilon_{m-M+i}+\frac{1}{M}\sum_{i=1}^M \left(f_-\l(\frac{m-M+i}{n}\r)-\hat{f}_{-}\l(\frac{m-M+i}{n}\r)\right)\right|\nonumber\\ 
    &\le \left|\frac{1}{M}\sum_{i=1}^M\varepsilon_{m-M+i}\right|\label{errorsum1}\\
    &\qquad +\left|\frac{1}{M}\sum_{i=1}^M \left(f_-\l(\frac{m-M+i}{n}\r)-\hat{f}_{-}\left(\frac{m-M+i}{n}\right)\right)\right|.\label{f-estiamtionerror1}
\end{align}
\end{subequations}
And for \(m_0 = \ceil{n\tau}+ 3N-1\), i.e.\ the test statistic after the true change point,
\begin{subequations}
\begin{align}\nonumber
    J_{m_0} &= \left|\frac{1}{M}\sum_{i=1}^M\varepsilon_{m_0-M+i} +\frac{1}{M}\sum_{i=1}^M \left(f_+\l(\frac{m_0-M+i}{n}\r)-\hat{f}_{-}\l(\frac{m_0-M+i}{n}\r)\right)\right|\\ \nonumber
    &= \left|\frac{1}{M}\sum_{i=1}^M\varepsilon_{m_0-M+i}\right.+\frac{1}{M}\sum_{i=1}^M \left(f_-\l(\frac{m_0-M+i}{n}\r)-\hat{f}_{-}\l(\frac{m_0-M+i}{n}\r)\right)\\ \nonumber
    &\qquad\qquad \left.+\frac{1}{ M}\sum_{i=1}^M \left(f_+\l(\frac{m_0-M+i}{n}\r)-f_-\l(\frac{m_0-M+i}{n}\r)\right)\right|\\ 
    &\ge \left|\frac{1}{M}\sum_{i=1}^M \left(f_+\l(\frac{m_0-M+i}{n}\r)-f_-\l(\frac{m_0-M+i}{n}\r)\right)\right|\label{f+-f-}\\
    &\qquad-\left|\frac{1}{M}\sum_{i=1}^M\varepsilon_{m_0-M+i}\right|\label{errorsum2}\\
    &\qquad\qquad -\left|\frac{1}{M}\sum_{i=1}^M \left(f_-\l(\frac{m_0-M+i}{n}\r)-\hat{f}_{-}\l(\frac{m_0-M+i}{n}\r)\right)\right|\label{f-estimationerror2}.
\end{align}
\end{subequations}
For the summand \eqref{f+-f-}, we have, because \(m_0-M- n\tau \le 3N-M\le N\),
\begin{align*}
    &\left|\frac{1}{M}\sum_{i=1}^M \left(f_+\l(\frac{m_0-M+i}{n}\r)-f_-\l(\frac{m_0-M+i}{n}\r)\right)\right| \\
    = & \left|\frac{1}{M}\sum_{i=1}^M \left(\l(\beta_+-\beta_-\r)\l(\frac{m_0-M+i}{n}-\tau\r)+\alpha_+-\alpha_-\right) \right|\\
    \ge & \l|\alpha_+-\alpha_-\r|-\l|\beta_+-\beta_-\r|\frac{1}{M}\sum_{i=1}^M\l(\frac{m_0-M+i}{n}-\tau\r)\\
    \ge &\l|\alpha_+-\alpha_-\r|-\l|\beta_+-\beta_-\r|\frac{1}{M}\sum_{i=1}^M\frac{N+i}{n} = \l|\alpha_+-\alpha_-\r|-\l|\beta_+-\beta_-\r| \frac{2N+M(M+1)}{2n} >c,
\end{align*}
for some \(c < |\alpha_+-\alpha_-| \) and $n$ large enough. The constant \(c\) exists, as \(0<\delta_0 \le  |\alpha_+-\alpha_-| \). 

We set 
\(Z_m\coloneqq \sum_{i=1}^M\varepsilon_{m-M+i}/\sqrt{M} \sim \mathcal{N}(0,1)\) and define the event
\begin{equation}\label{Ajump}
    \mathcal{A}\coloneqq \left\{\max_{k< m\le n} |Z_m| < \sqrt{10 \log n},\; \max_{1\le i\le n} \l|f_-\l(\frac{i}{n}\r)-\hat f_{-}\l(\frac{i}{n}\r)\r|< \frac{c}{10}\right\}. 
\end{equation}
For the summands \eqref{errorsum1} and \eqref{errorsum2} under the event \(\mathcal{A}\) with \(b = 10^3/(2c^2)\) it holds
\begin{equation}\label{jumperrorbound}
    \left|\frac{1}{M}\sum_{i=1}^M\varepsilon_{m-M+i}\right|\le
    \max_{k< m \le n} \left|\frac{Z_m}{\sqrt{M}}\right| < \frac{\sqrt{10 \log n}}{\sqrt{M}}\le \sqrt{\frac{10\log n}{2N}} = \sqrt{\frac{10\log n}{2 b \log n}} = \sqrt{\frac{10 c^2}{10^3}}= \frac{c}{10}.
\end{equation}
For the parts \eqref{f-estiamtionerror1} and \eqref{f-estimationerror2} under the event \(\mathcal{A}\) it holds
\begin{align}
    &\left|\frac{1}{M}\sum_{i=1}^M \left(f_-\l(\frac{m-M+i}{n}\r)-\hat{f}_{-}\l(\frac{m-M+i}{n}\r)\right)\right|\nonumber \\
    &\qquad\le \frac{1}{M}\sum_{i=1}^M\left|f_-\l(\frac{m-M+i}{n}\r)-\hat{f}_{-}\l(\frac{m-M+i}{n}\r)\right|<\frac{c}{10}.\label{f-estimationerrorbound}
\end{align}
Let \(c < \delta_0\). Under the event \(\mathcal{A}\), for \(m < n\tau\), it holds that  \(J_m < c/10 + c/10 < 4c/5 = \rho_J\) and \(J_{m_0} > c - c/10 - c/10 = 4c/5=\rho_J\), and thus, \(n\tau \le n\hat\tau_{J,n} \le m_0\).

For \(Z\sim\mathcal{N}(0,1)\), using Mill's ratio, we obtain
\begin{equation}\label{millsratio}\mathbb{P}\left(|Z|\ge y\right) = 2 \mathbb{P}\left( Z \ge y\right)\le\exp\left\{-\frac{y^2}{2}\right\},
\end{equation}
and therefore using the union bound
\begin{equation}
    \mathbb{P}_{\tau}\left(\max_{k< m \le n} |Z_m|\ge y\right) \le n\exp\left\{-\frac{y^2}{2}\right\}.\label{pmaxz}
\end{equation}
When we set \(y = \sqrt{10 \log n}\), we have 
\begin{equation}\label{boundP(A_1)}
    \mathbb{P}_{\tau}\left(\max_{k< m \le n} |Z_m|\ge \sqrt{10 \log n}\right)\le n\exp\left\{-5\log n\right\}= n^{-4} .
\end{equation}

By \citet[Theorem~7.5]{korostelev}, we have
\[\hat{\alpha}\sim \mathcal{N}\l(\alpha_--\beta_-\tau,\frac{4k+2}{k^2-k}\r)\text{ and }\hat{\beta}\sim \mathcal{N}\l(\beta_-,\frac{12n^2}{k^3-k}\r).\] 
Define 
\[Z_\alpha = \l(\alpha_--\beta_-\tau-\hat{\alpha}\r)\sqrt{\frac{k^2-k}{4k+2}}\sim \mathcal{N}\l(0,1\r)\text{ and }Z_\beta = \l(\beta_--\hat{\beta}\r)\sqrt{\frac{k^3-k}{12n^2}}\sim \mathcal{N}\l(0,1\r).\]
Now we can use the inequality \eqref{millsratio} to obtain for any \(\epsilon > 0\)
\begin{align}
   &\P\l(\max_{1\le i \le n}\l|f_-\l(\frac{i}{n}\r)-\hat{f}_-\l(\frac{i}{n}\r)\r|\ge\epsilon\r) = \P\l(\l\{\l|f_-(0)-\hat{f}_-(0)\r|\ge \epsilon\r\}\cup\l\{ \l|f_-(1)-\hat{f}_-(1)\r|\ge \epsilon\r\}\r)\nonumber\\
    =\, &\P\l(\l\{\l|\alpha_--\beta_-\tau-\hat{\alpha}\r|\ge \epsilon\r\} \cup \l\{\l|\alpha_-+\beta_-\l(1-\tau\r) -\hat{\alpha}-\hat{\beta}\r|\ge\epsilon\r\}\r)\nonumber\\
    \le\,&\P\l(\l\{\l|\alpha_--\beta_-\tau-\hat{\alpha}\r|\ge\frac{\epsilon}{2}\r\} \cup\l\{ \l|\beta_--\hat{\beta}\r|\ge\frac{\epsilon}{2}\r\}\r)\nonumber\\
    \le\,&\P\l(\l|\alpha_--\beta_-\tau-\hat{\alpha}\r|\ge \frac{\epsilon}{2}\r)+\P\l( \l|\beta_--\hat{\beta}\r|\ge\frac{\epsilon}{2}\r)\nonumber\\
    = \,&\P\l(\l|Z_\alpha\r|\ge \frac{\epsilon}{2}\sqrt{\frac{k^2-k}{4k+2}}\r)+\P\l(\l|Z_\beta\r|\ge\frac{\epsilon}{2}\sqrt{\frac{k^3-k}{12n^2}}\r)\nonumber\\
    \label{estimationerrorbound}
    \le\, & \exp \l\{-\frac{\epsilon^2\l(k^2-k\r)}{32k+16}\r\}+\exp\l\{-\frac{\epsilon^2(k^3-k)}{96 n^2}\r\},
\end{align}
and therefore 
\begin{align}
    \P\l(\max_{1\le i \le n}\l|f_-\l(\frac{i}{n}\r)-\hat{f}_-\l(\frac{i}{n}\r)\r|\ge\frac{c}{10}\r) &\le \exp \l\{-\frac{c^2\l(c^2n^2-c n\r)}{3200c n+16}\r\}+\exp\l\{-\frac{c^2(c^3n^2-c)}{9600 n}\r\}\nonumber\\
    &\ll n^{-4}.\label{boundP(A_2)}
\end{align}

We can combine \eqref{boundP(A_1)} and \eqref{boundP(A_2)} to get
\begin{equation*}
    \P\l(\mathcal{A}^c\r)\le \P\l(\max_{1\le m\le M}|Z_m|\ge y\r)+ \P\l(\max_{1\le i \le n}\l|f_-\l(\frac{i}{n}\r)-\hat{f_-}\l(\frac{i}{n}\r)\r|\ge \frac{c}{10}\r) \le n^{-4}.
\end{equation*}
To summarize, we have shown that under the event \(\mathcal{A}\), it holds that \(0\le\hat \tau_{J,n}-\tau\le3N/n\) and  \(\P\l(\mathcal{A}^c\r)\le n^{-4}\). Thus, since \(1/\l(n^2 \log^2 n\r)< 1\) for \(n\ge2\),
\begin{align*}&\max_{\theta_0\in\Theta^J}\mathbb{E}_{\tau}\left[\left(\frac{n}{\log n}\l(\hat\tau_{J,n}-\tau\r)\right)^2\right]\\ &\qquad = \max_{\theta_0\in\Theta^J}\left(\mathbb{E}_{\tau}\left[\left(\frac{n}{\log n}\l(\hat\tau_{J,n}-\tau\r)\right)^2\mathbb{I}(\mathcal{A})\right]+ \mathbb{E}_{\tau}\left[\left(\frac{n}{\log n}\l(\hat\tau_{J,n}-\tau\r)\right)^2\mathbb{I}(\mathcal{A}^c)\right]\right)\\
    &\qquad \le \max_{\theta_0\in\Theta^J}\left(\mathbb{E}_{\tau}\left[\left(\frac{3N}{\log n}\right)^2\mathbb{I}(\mathcal{A})\right]+ \mathbb{E}_{\tau}\left[\left(\frac{n}{\log n}\right)^2\mathbb{I}(\mathcal{A}^c)\right]\right)\\
    &\qquad \le \left(\frac{3N}{\log n}\right)^2+\left(\frac{n}{\log n}\right)^2 n^{-4}\le 9b^2 + 1 = 9\cdot\frac{10^6}{4c^4} +1.
\end{align*}
Therefore, the statement about a jump is true with \(r_J^* = 9(10^6/4\delta_0^4) +1. \)

{\rev We additionally note that the choice $k \asymp n$ is required to establish  \eqref{boundP(A_2)}. However, by leveraging~\eqref{estimationerrorbound}, this requirement can be relaxed to  $k \asymp n^{2/3}\log(n)^{1/3}$.}

\paragraph{For a kink} We assume that \(\l|\beta_+-\beta_-\r|\ge\delta_0\) and \(\l|\alpha_+-\alpha_-\r|<\delta_0\). Note that \(N \coloneqq bn^{2/3}\log^{1/3} n \), with \(b\) positive and independent of \(n\),  and \(M \coloneqq 2N + (m\mod N)\). The threshold is \(\rho_K = 4c /(5 n)\) for some fixed \(c >0\).
Let \(d_M \coloneqq M(M+1)(2M+1)/6\).
For \(m < n\tau\), it holds
\begin{subequations}
    \begin{align}\nonumber
    K_m &= \left|\frac{1}{d_M}\sum_{i=1}^Mi\varepsilon_{m-M+i}+\frac{1}{d_M}\sum_{i=1}^M i\left(f_-\l(\frac{m-M+i}{n}\r)-\hat{f}_{-}\l(\frac{m-M+i}{n}\r)\right)\right|\\ 
    &\le \left|\frac{1}{d_M}\sum_{i=1}^Mi\varepsilon_{m-M+i}\right|\label{kinkerrorsum1}\\
    &\qquad +\left|\frac{1}{d_M}\sum_{i=1}^M i\left(f_-\l(\frac{m-M+i}{n}\r)-\hat{f}_{-}\l(\frac{m-M+i}{n}\r)\right)\right|.\label{kinkf-estiamtionerror1}
\end{align}
\end{subequations}
And for \(m_0 \coloneqq \ceil{n\tau}+3N-1\), we get
\begin{subequations}
    \begin{align}\nonumber
    K_{m_0} &= \left|\frac{1}{d_M}\sum_{i=1}^Mi\varepsilon_{m_0-M+i}+\frac{1}{d_M}\sum_{i=1}^M i\left(f_+\l(\frac{m_0-M+i}{n}\r)-\hat{f}_{-}\l(\frac{m_0-M+i}{n}\r)\right)\right|\\ \nonumber
    &= \l|\frac{1}{d_M}\sum_{i=1}^Mi\varepsilon_{m_0-M+i}+\frac{1}{d_M}\sum_{i=1}^M i\left(f_-\l(\frac{m_0-M+i}{n}\r)-\hat{f}_{-}\l(\frac{m_0-M+i}{n}\r)\right)\r.\\ \nonumber
    &\qquad \left.+\frac{1}{ d_M}\sum_{i=1}^M i\left(f_+\l(\frac{m_0-M+i}{n}\r)-f_-\l(\frac{m_0-M+i}{n}\r)\right)\right|\\ 
    &\ge \left|\frac{1}{d_M}\sum_{i=1}^M i\left(f_+\l(\frac{m_0-M+i}{n}\r)-f_-\l(\frac{m_0-M+i}{n}\r)\right)\right|\label{kinkf+-f-}\\
    &\qquad-\left|\frac{1}{d_M}\sum_{i=1}^M i\varepsilon_{m_0-M+i}\right|\label{kinkerrorsum2}\\
    &\qquad\qquad -\left|\frac{1}{ d_M}\sum_{i=1}^M i\left(f_-\l(\frac{m_0-M+i}{n}\r)-\hat{f}_{-}\l(\frac{(m_0-M+i}{n}\r)\right)\right|.\label{kinkf-estimationerror2}
\end{align}
\end{subequations}

As there is a kink and \(m_0-M-n\tau \ge 0\). we get for the summand \eqref{kinkf+-f-}, in the case of \(\l|\alpha_+-\alpha_-\r| = 0\),
\begin{align*}
    &\left|\frac{1}{d_M}\sum_{i=1}^M i\left(f_+\l(\frac{m_0-M+i}{n}\r)-f_-\l(\frac{m_0-M+i}{n}\r)\right)\right|\\
    &\qquad= \l| \frac{1}{d_M}\sum_{i=1}^Mi\l(\l(\beta_+-\beta_-\r) \l(\frac{m_0-M+i}{n}-\tau\r)\r)\r|\\
    &\qquad \ge \frac{\l|\beta_+-\beta_-\r|}{d_M}\sum_{i=1}^Mi\l(\frac{m_0-M+i}{n}-\tau\r)\ge \frac{\l|\beta_+-\beta_-\r|}{d_Mn}\sum_{i=1}^M i^2 > \frac{c}{n},
 \end{align*}
for \(\l|\beta_+-\beta_-\r| > c\).
In the case of \(\l|\alpha_+-\alpha_-\r|> 0\), as \(m_0-M-n\tau \le N\),  the summand \eqref{kinkf+-f-} can be bounded as follows
\begin{align*}
    &\left|\frac{1}{d_M}\sum_{i=1}^M i\left(f_+\l(\frac{m_0-M+i}{n}\r)-f_-\l(\frac{m_0-M+i}{n}\r)\right)\right|\\
= \,&\l| \frac{1}{d_M}\sum_{i=1}^Mi\l(\l(\beta_+-\beta_-\r) \l(\frac{m_0-M+i}{n}-\tau\r)+\alpha_+-\alpha_-\r)\r|\\
 \ge\,& \frac{\l|\alpha_+-\alpha_-\r|}{d_M}\sum_{i=1}^M i - \frac{\l|\beta_+-\beta_-\r|}{d_M}\sum_{i=1}^M i\l(\frac{m_0-M+i}{n}-\tau\r)\\
     \ge\,& \frac{3\l|\alpha_+-\alpha_-\r|}{2M+1}-\frac{\l|\beta_+-\beta_-\r|}{d_Mn}\sum_{i=1}^Mi(i+N) = \frac{3\l|\alpha_+-\alpha_-\r|}{2M+1} - \frac{\l|\beta_+-\beta_-\r|}{n}\l(1+\frac{3N}{2M+1}\r)
     > \frac{c}{n},
\end{align*}
for \(n\) large enough. 
Set \(Z_m \coloneqq (1/\sqrt{d_M})\sum_{i=1}^M i \varepsilon_{m-M+i} \) with \(\varepsilon_i\) independently standard normal distributed. Next we define the event 
\begin{equation}\label{Akink}
    \mathcal{A}\coloneqq \left\{\max_{k< m\le n} |Z_m| < \sqrt{8 \log n}, \max_{1\le i\le n} \l|f_-\l(\frac{i}{n}\r)-\hat f_{-}\l(\frac{i}{n}\r)\r|< \frac{c\l(2M+1\r)}{30n}\right\}.
\end{equation}
For the parts \eqref{kinkerrorsum1} and \eqref{kinkerrorsum2}, it  holds under the event \(\mathcal{A}\) with \(b = \l(300/c^2\r)^{1/3}\) that
\begin{align}
    \left|\frac{1}{d_M}\sum_{i=1}^Mi\varepsilon_{m-M+i}\right| &\le \max_{k< m\le n} \frac{1}{\sqrt{d_M}}\l|Z_m\r| < \frac{\sqrt{8 \log n}}{\sqrt{d_M}} = \sqrt{\frac{48\log n}{2M^3+3M^2+M}}\le \sqrt{\frac{24\log n}{(2N)^3}}\nonumber\\
    &= \sqrt{\frac{3\log n}{b^3n^2\log n}} = \frac{c}{10n}.\label{kinkerrorsumbound}
\end{align}
And for the parts \eqref{kinkf-estiamtionerror1} and \eqref{kinkf-estimationerror2}, it holds under the event \(\mathcal{A}\) that
\begin{multline}
    \left|\frac{1}{ d_M}\sum_{i=1}^M i\left(f_-\l(\frac{m-M+i}{n}\r)-\hat{f}_{-}\l(\frac{m-M+i}{n}\r)\right)\right| < \frac{1}{ d_M}\sum_{i=1}^M i\frac{c\l(2M+1\r)}{30n}\\
    = \frac{6}{M\l(M+1\r)\l(2M+1\r)}\frac{c\l(2M+1\r)}{30n}\frac{M\l(M+1\r)}{2} = \frac{c}{10n}.\label{kinkf-estimationerrorbound}
\end{multline}
Thus, it holds 
\[\mathcal{A}\subset \l\{K_m <\frac{c}{5n} \text{ for } m<n\tau\text{ and } K_m \ge \frac{4c}{5n} \text{ for } m=\ceil{n\tau}+3N-1 \r\}\subset \l\{n\tau\le n \hat\tau_{K,n}\le n\tau+3N\r\}\] and therefore under the event \(\mathcal{A}\) 
\begin{equation}\label{kinktau-theta}
    0\le\hat\tau_{K,n}-\tau\le \frac{3N}{n}.
\end{equation}
Now we can use \eqref{pmaxz} to get
\begin{equation*}
    \P_{\tau}\l(\max_{k< m\le n}|Z_m|\ge \sqrt{8\log n}\r) \le n \exp\l\{-4\log n\r\} = n^{-3}.
\end{equation*}
With the inequality \eqref{estimationerrorbound} it follows that
\begin{align}
        &\P_{\tau}\l(\max_{1\le i \le n}\l|f_-\l(\frac{i}{n}\r)-\hat{f}_-\l(\frac{i}{n}\r)\r|\ge\frac{c\l(2M+1\r)}{30n}\r)\nonumber\\
 \le\,& \exp \l\{-\frac{c^2(2M+1)^2\l(c^2n^2-c n\r)}{300n^2(8c n+4)}\r\}+\exp\l\{-\frac{c^2(2M+1)^2(c^3n^2-c)}{28800 n^3}\r\},\nonumber\\
         \asymp\,& \exp\l\{-n^{1/3}\log^{2/3} n\r\} \ll n^{-3}.\label{eq:fmhat}
\end{align}
Thus, it holds
\begin{align*}
\P_{\tau}\l(\mathcal{A}^c\r)\le \P_{\tau}\l(\max_{1\le m\le M}|Z_m|\ge \sqrt{8\log M}\r)  +\P_{\tau}\l(\max_{1\le i \le n}\l|f_-\l(\frac{i}{n}\r)-\hat{f}_-\l(\frac{i}{n}\r)\r|\ge\frac{c\l(2N+1\r)}{10n}\r)\le n^{-3}.
\end{align*}
Based on this and \eqref{kinktau-theta}, we get
\begin{align*}&\max_{\theta_0\in\Theta^K}\mathbb{E}_{\tau}\left[\left(\frac{n^{1/3}}{\log^{1/3}n}\l(\tau_{K,n}-\tau\r)\right)^2\right]\\ 
=\,& \max_{\theta_0\in\Theta^K}\left(\mathbb{E}_{\tau}\left[\left(\frac{n^{1/3}}{\log^{1/3}n}\l(\tau_{K,n}-\tau\r)\right)^2\mathbb{I}(\mathcal{A})\right]+ \mathbb{E}_{\tau}\left[\left(\frac{n^{1/3}}{\log^{1/3}n}\l(\tau_{K,n}-\tau\r)\right)^2\mathbb{I}(\mathcal{A}^c)\right]\right)\\
 \le \,&\max_{\theta_0\in\Theta}\left(\mathbb{E}_{\tau}\left[\left(\frac{3N}{n^{2/3}\log^{1/3}n}\right)^2\mathbb{I}(\mathcal{A})\right]+ \mathbb{E}_{\tau}\left[\left(\frac{n^{1/3}}{\log^{1/3}n}\right)^2\mathbb{I}(\mathcal{A}^c)\right]\right)\\
 \le\,& \left(\frac{3N}{n^{2/3}\log^{1/3}n}\right)^2+\left(\frac{n^{1/3}}{\log^{1/3}n}\right)^2 n^{-3} \le 9b^2 + 1 = 9 \cdot \left(\frac{300}{c^2}\right)^{2/3} + 1,
\end{align*}
by the fact \(1/\log^{2/3}n\le 1\) for \(n>2\).
Thus, the statement for a kink holds with \(r^*_K = 9\left(300/\delta_0^2\right)^{2/3} + 1 \).

{\rev
Analogous to the jump case, the requirement on the size of historical data can be relaxed from $k \asymp n$
 to $k \asymp n^{8/9}\log(n)^{1/9}$; thanks to~\eqref{estimationerrorbound}, this relaxed choice still guarantees~\eqref{eq:fmhat}.
}
\end{proof}

{\rev
\begin{proof}[Proof of \Cref{theorem:falsedetectiontypebound}]
Note that \(N_J\coloneqq b_J \log n\), \(M_J\coloneqq 2N_J + (m \mod N_J)\), \(b_J \coloneqq 10^3/ \l(2c^2\r)\) and \(\rho_J \coloneqq 4c/5\) in the jump case, while \(N_K\coloneqq b_k n^{2/3}\log^{1/3} n\), \(M_K\coloneqq 2N_K + (m \mod N_K)\), \(b_K \coloneqq (300/c^2)^{1/3}\) and \(\rho_K = 4c/(5n)\) in the kink case.
We set \(Z_{J,m}\coloneqq \sum_{i=1}^M\varepsilon_{m-M+i}/\sqrt{M} \sim \mathcal{N}(0,1)\) and \(Z_{k,m} \coloneqq (1/\sqrt{d_M})\sum_{i=1}^M i \varepsilon_{m-M+i} \).
We define events \(\mathcal{A}_J = \mathcal{A}\) in \eqref{Ajump} and  \(\mathcal{A}_K = \mathcal{A}\) in \eqref{Akink}.
In the proof of \Cref{theorem:upperbound} we showed that \(\P\l(\mathcal{A}^c_J\r)\le n^{-4}\) and \(\P\l(\mathcal{A}^c_K\r)\le n^{-3}\), and thus \(\P\l(\l(\mathcal{A}_J\cap \mathcal{A}_K\r)^c\r)\ll n^{-3}\).
Also the corresponding detector will not make a false detection under the events \(\mathcal{A}_J\) and \(\mathcal{A}_K\).

\paragraph{For a jump}
Fix \(\theta_0 \in \Theta_{\delta_0}^J\). We only need to show \(\{\hat\tau_{K,n}<\hat\tau_{J,n}\}\subseteq \{\mathcal{A}_J\cap \mathcal{A}_K\}^c\). Note that we already showed in the proof of \Cref{theorem:upperbound} that under a jump and the event \(\mathcal{A}_J\) is \(\hat\tau_{J,n}\le \tau_0 + 3N_J/n \) and that under the event \(\mathcal{A}_K\) is \(\tau_0\le \hat\tau_{K,n}\)
Now we choose \(m\in \l[\tau_0n,\tau_0n+3N_J\r]\):
\begin{subequations}
    \begin{align}\nonumber
    K_{m} &= \left|\frac{1}{d_{M_K}}\sum_{i=1}^{M_K}i\varepsilon_{m-M_K+i}+\frac{1}{d_{M_K}}\sum_{i=1}^{M_K} i\left(f_+\l(\frac{m-M_K+i}{n}\r)-\hat{f}_{-}\l(\frac{m-M_K+i}{n}\r)\right)\right|\nonumber\\ 
    &\le \left|\frac{1}{d_{M_K}}\sum_{i=1}^{M_K} i\left(f_+\l(\frac{m-M_K+i}{n}\r)-f_-\l(\frac{m-M_K+i}{n}\r)\r)\mathbb{I}\l(i> \tau_0 n-m + M_K\right)\right|\label{kinkf+-f-3}\\
    &\qquad+\left|\frac{1}{d_{M_K}}\sum_{i=1}^{M_K} i\varepsilon_{m-M_K+i}\right|\label{kinkerrorsum3}\\
    &\qquad\qquad +\left|\frac{1}{ d_{M_K}}\sum_{i=1}^{M_K} i\left(f_-\l(\frac{m-M_K+i}{n}\r)-\hat{f}_{-}\l(\frac{(m-M_K+i}{n}\r)\right)\right|,\label{kinkf-estimationerror3}
\end{align}
\end{subequations}
where both terms in  \eqref{kinkerrorsum3} and  \eqref{kinkf-estimationerror3} can be bounded by \(c/(10n)\), see \eqref{kinkerrorsumbound} and  \eqref{kinkf-estimationerrorbound}. For \eqref{kinkf+-f-3} we get
\begin{align*}
    &\left|\frac{1}{d_{M_K}}\sum_{i=1}^{M_K} i\left(f_+\l(\frac{m-M_K+i}{n}\r)-f_-\l(\frac{m-M_K+i}{n}\r)\r)\mathbb{I}\l(i> \tau_0 n-m + M_K\right)\right|\\
=\,&\l|\frac{1}{d_{M_k}} \sum_{i=\tau_0 n-m+M_K}^{M_K} i\l( \l(\beta_+ - \beta_-\r) \l(\frac{m-M_K+i}{n} - \tau_0\r) + \alpha_+ -\alpha_- \r)\r|\\
 \le \,&\frac{1}{d_{M_K}} \sum_{i=0}^{m-\tau_0n} \left(i^2 \frac{\l|\beta_+-\beta_-\r|}{n} + i\l|\alpha_+ - \alpha_-\r|\right) \le \frac{1}{d_{M_K}} \sum_{i=0}^{3N_J}\left( i^2 \frac{\l|\beta_+-\beta_-\r|}{n} + i\l|\alpha_+ - \alpha_-\r|\right)\\
     \le \,&\frac{3N_J(3N_J+1)(6N_J+1)}{3N_K(3N_K+1)(6N_K+1)}\frac{\l|\beta_+-\beta_-\r|}{n} + \frac{9N_J(N_J+1)}{3N_K(3N_K+1)(6N_K+1)} \l|\alpha_+-\alpha_-\r|\\
\asymp\,& \frac{\log^2n}{n^3} \l|\beta_+ - \beta_-\r| + \frac{\log n}{n^2} \l|\alpha_+ -\alpha_-\r|.
\end{align*}

Now we can bound \eqref{kinkf+-f-3} by \({c}/{(5n)}\) for \(n\) large enough and finally get that \(K_m \le c/(5n)+ c/(10n) + c/(10n) = 2/(5n)<\rho_K\). Hence we get that under both events \(\mathcal{A}_J\) and \(\mathcal{A}_K\) is \(\hat\tau_{k,n} > \tau_0  + 3N_J/n \ge \hat\tau_{J,n}\).

\paragraph{For a kink}

For \(\theta_0 \in \Theta^K_{\delta_0}\), we want to show \(\{\hat\tau_{J,n}<\hat\tau_{K,n}\}\subseteq \{\mathcal{A}_J\cap \mathcal{A}_K\}^c\). In the proof of \Cref{theorem:upperbound} we showed that  \(\hat\tau_{K,n}\le\tau_0+3N_K/n\)  in this case  under the event \(\mathcal{A}_K\). For \(m\in \l[\tau_0n,\tau_0 n+3N_K\r]\),
\begin{subequations}
\begin{align}\nonumber
    J_{m} &= \left|\frac{1}{M_J}\sum_{i=1}^{M_J}\varepsilon_{m-M_J+i} +\frac{1}{M_J}\sum_{i=1}^{M_J} \left(f_+\l(\frac{m-M_J+i}{n}\r)-\hat{f}_{-}\l(\frac{m-M_J+i}{n}\r)\right)\right|\\
    &\le \left|\frac{1}{M_J}\sum_{i=1}^{M_J} \left(f_+\l(\frac{m-M_J+i}{n}\r)-f_-\l(\frac{m-M_J+i}{n}\r)\right)\mathbb{I}\l(i> \tau_0 n-m + M_J\r)\right|\label{f+-f-3}\\
    &\qquad+\left|\frac{1}{M_J}\sum_{i=1}^{M_J}\varepsilon_{m-M_J+i}\right|\label{errorsum3}\\
    &\qquad\qquad +\left|\frac{1}{M_J}\sum_{i=1}^{M_J} \left(f_-\l(\frac{m-M_J+i}{n}\r)-\hat{f}_{-}\l(\frac{m-M_J+i}{n}\r)\right)\right|\label{f-estimationerror3}.
\end{align}
\end{subequations}

We have shown in \eqref{jumperrorbound} that under the event \(\mathcal{A}_J\) \eqref{errorsum3} is bounded by \(c/10\) and in \eqref{f-estimationerrorbound} that \eqref{f-estimationerror3} is bounded by \(c/10\). For the part \eqref{f+-f-3} we get
\begin{align*}
    &\left|\frac{1}{M_J}\sum_{i=1}^{M_J} \left(f_+\l(\frac{m-M_J+i}{n}\r)-f_-\l(\frac{m-M_J+i}{n}\r)\right)\mathbb{I}\l(i> \tau_0 n-m + M_J\r)\right|\\
    \le\,& \frac{1}{M_J} \sum_{i=1}^{M_J} \l| \beta_+ - \beta_-\r| \l(\frac{m-M_J+i}{n}-\tau_0\r)+\l|\alpha_+-\alpha_-\r| \\
     \le \,&\frac{1}{M_J} \sum_{i=1}^{M_J} \l|\beta_+ - \beta_-\r| \l(\frac{3N_K-M_J+i}{n}\r) + \l|\alpha_+-\alpha_-\r| \le \l|\beta_+ - \beta_-\r| \frac{3N_K}{n} + \l|\alpha_+ -\alpha_-\r|\\
    =\,&  \l|\beta_+ - \beta_-\r| \frac{3b_k\log^{1/3}n}{n^{1/3}}  + \l|\alpha_+ -\alpha_-\r| < \frac{c}{5},
\end{align*}
for \(n\) large enough, since we have that \(\l|\alpha_+-\alpha_-\r| < \delta_0 \).
As a result we get \(J_m < c/10 + c/10 + c/5 = 2c/5 < \rho_J\) and now have under the event \(\l\{\mathcal{A}_J \cap \mathcal{A}_K\r\}\) that \(\hat\tau_{J,n}>\tau_0 + 3N_K/n \ge \hat \tau_{K,n}\).
\end{proof}

\begin{proof}[Proof of \Cref{corollary:detection delay}]
We prove the bound on the expected detection delay by utilizing \eqref{eq:typeIandtypeIIerror} and \eqref{eq:detdelbound}. By these equations we get for any detector \(\tau_n\) the inequality
\begin{equation}\label{eq:detdelbound2}
\bigl(\E_\tau\l[\hat\tau_n-\tau\mid\hat\tau_n\ge\tau\r]\bigr)^2 \le \frac{\E_\tau\l[\l(\hat\tau_n-\tau\r)^2\r]}{\bigl(1 - \P_\tau(\hat\tau_n< \tau)\bigr)}.
\end{equation}

In the proof of \Cref{theorem:upperbound} we showed that \(\P_\tau\l(\hat\tau_{J,n}<\tau\r)\le n^{-3} \) and
\[\sup_{\theta = (\tau,\alpha_{-},\alpha_{+},\beta_{-},\beta_{+})\in\Theta^J_{\delta_0}}\mathbb{E}_{\tau}\left[\left(\frac{n}{\log n}\l(\hat\tau_{J,n}-\tau\r)\right)^2\right]\;\le\; r^*_J < \infty.\]
By plugging this into \eqref{eq:detdelbound2} we get for \(n\) large enough a bound on the expected detection delay
\[\sup_{\theta = (\tau,\alpha_{-},\alpha_{+},\beta_{-},\beta_{+})\in\Theta^J_{\delta_0}}\frac{n}{\log n} \E_\tau\l[\hat\tau_{J,n}-\tau\mid\hat\tau_{J_n}\ge\tau\r]\le \sqrt{r_J^*}.\] 
To bound the average run length, we recall that  \(\hat\tau_{J,n}\) is set to \(1\) if no detection is made. Thus the bound on the average run length under the null is 
\[\E_0\l[\hat\tau_{J,n}\r]\ge \P_0(\hat\tau_{J,n} = 1)\ge 1 - n^{-3}.\]
The proof for the kink case is analogous. We have already shown that \(\P\l(\hat\tau_{K,n}\le\tau\r)\le n^{-4} \) and 
\[
        \sup_{\theta = (\tau,\alpha_{-},\alpha_{+},\beta_{-},\beta_{+})\in\Theta^K_{\delta_0}}\mathbb{E}_{\tau}\left[\left(\frac{n^{1/3}}{\log^{1/3}n}\l(\hat\tau_{K,n}-\tau\r)\right)^2\right]\le r^*_K < \infty.
        \]
With these results we bound the expected detection delay by
\[\sup_{\theta = (\tau,\alpha_{-},\alpha_{+},\beta_{-},\beta_{+})\in\Theta^K_{\delta_0}} \frac{n^{1/3}}{\log^{1/3} n} \E_\tau\l[\hat\tau_{K,n}-\tau|\hat\tau_{K,n}\ge\tau\r]\le \sqrt{r_K^*},\]
and the average run length by
\[\E_0\l[\hat\tau_{K,n}\r]\ge 1 - n^{-4}.\]
Combining the above statements shows the claim for the FLOC detector $\hat\tau_n = \min\{\hat\tau_{J,n}, \hat\tau_{K,n}\}$.
\end{proof}
}

\begin{proof}[Proof of \Cref{theorem:lowerbound}]
    We are going to prove the statement by contradiction. The proof follows a similar structure to that of Theorem 6.12 in \cite{korostelev}, which establishes the result for the special  case of a jump and \(\beta_- = \beta_+ = 0\).
    The initial steps are identical for the cases of a kink and a jump. Let \(\psi_n\) be a sequence satisfying \(\frac{1}{n}\le \psi_n \ll 1\).

    We fix \(\alpha_- = \beta_- = 0\) and some \(\alpha_+\) and \(\beta_+\), this is possible, since a lower rate obtained over the fixed parameters is also a lower rate for the whole parameter space. 
    We set \(M \coloneqq \floor{\l(1-2\delta_0\r)/(3 b \psi_n)}\), with \(b\) a positive constant independent of \(n\), to be determined later.  
    We choose \(t_0,\ldots, t_M\in (\delta_0,1-\delta_0)\) such that \(t_j-t_{j-1}= 3 b \psi_n\) for \(j = 0,\ldots,M\) and \(M>1\). 
    As \(\frac{1}{n}\le \psi_n \ll 1\), this is possible. 
    Now we assume that
    \[\liminf_{n\to\infty}\inf_{\hat{\tau}_n\in \mathcal{T}}\max_{\tau_0\in(\delta_0,1-\delta_0)} \mathbb{E}_{\tau_0}\left[\left(\frac{\hat{\tau}_n-\tau_0}{\psi_n}\right)^2\right]=0.\] 
    Then, there exists a family of Markov stopping times \((\tilde{\tau}_n)_{n\in\mathbb{N}}\) such that 
    \begin{equation}\label{lowerassumption}
        \lim_{n\to\infty} \max_{0\le j \le M}\mathbb{E}_{t_j}\left[\left(\frac{\tilde{\tau}_n-t_j}{\psi_n}\right)^2\right] = 0.
    \end{equation}
    By Markov's inequality, we have
    \[\mathbb{P}_{t_j}\left(\left|\tilde{\tau}_n- t_j\right|> b \psi_n\right)\le \frac{\mathbb{E}_{t_j}\left[\left|\tilde{\tau}_n-t_j\right|^2\right]}{(b\psi_n)^2},\]
    and further, by \eqref{lowerassumption},
    \[\lim_{n\to\infty}\max_{0\le j\le M} \mathbb{P}_{t_j}\left(\left|\tilde{\tau}_n-t_j\right|>b\psi_n\right) = 0.\]
    Hence, for \(n\) large enough,
    \begin{equation}\label{P(Bj)}
        \max_{0\le j\le M} \mathbb{P}_{t_j}\left(\left|\tilde{\tau}_n-t_j\right|>b\psi_n\right)\le \frac{1}{4}.
    \end{equation}
    Let \(\mathcal{B}_j\coloneqq \left\{\left|\tilde{\tau}_n-t_j\right|\le b \psi_n\right\}\).
    The distance between any two points of \(t_0,\ldots,t_M\) is greater than \(3b\psi_n\), so the events  \(\mathcal{B}_j\) are all mutually exclusive.
    Note also that \(\bigcup_{j=0}^{M-1}\mathcal{B}_j\subset \left\{\left|\tilde{\tau} -t_M\right| > b\psi_n\right\}\). We denote
    \[\mathbf{d}\mathbb{P}_{t_j} = \frac{1}{\left(\sqrt{2\pi}\right)^n} \prod_{i=1}^n \exp\left\{-\frac{1}{2}\left(X_i-f_{t_j}\right)^2\right\},\]
    the likelihood function of \(X_1, \ldots,X_n\) under the change point \(t_j\). Next we can calculate
    \begin{align*}
        &\sum_{j=0}^{M-1}\mathbb{E}_{t_j}\left[\frac{\mathbf{d}\mathbb{P}_{t_M}}{\mathbf{d}\mathbb{P}_{t_j}} \mathbb{I}\left(\mathcal{B}_j\right)\right] =\sum_{j=0}^{M-1} \mathbb{E}_{t_M}\left[ \mathbb{I}\left(\mathcal{B}_j\right)\right]= \sum_{j=0}^{M-1}\mathbb{P}_{t_M}\left(\mathcal{B}_j\right)= \mathbb{P}_{t_M}\left(\bigcup_{j=0}^{M-1}\mathcal{B}_j\right)\\
        \le\,& \mathbb{P}_{t_M}\left(\left|\tilde{\tau}_n-t_M\right| > b\psi_n\right) \le \max_{0\le j\le M} \mathbb{P}_{t_j}\left(\left|\tilde{\tau}_n-t_j\right|>b\psi_n\right).
    \end{align*}
    And by \eqref{P(Bj)}, we have, for \(n\) large enough,
    \begin{equation}
    \label{one}
        \sum_{j=0}^{M-1}\mathbb{E}_{t_j}\left[\frac{\mathbf{d}\mathbb{P}_{t_M}}{\mathbf{d}\mathbb{P}_{t_j}} \mathbb{I}\left(\mathcal{B}_j\right)\right]\le \frac{1}{4}.
    \end{equation}
    Now we can calculate the likelihood ratio.
    For notational convenience, we write \(f_j(i) \coloneqq f_{t_j}\l(i/n\r)\) and \(f_M(i) \coloneqq f_{t_M}\l({i}/{n}\r) \).
    Note that \(f_j(i) = 0 = f_M(i)\) for \(i\le t_j<t_M\). Then
    \begin{align*}
       \frac{\mathbf{d}\mathbb{P}_{t_M}}{\mathbf{d}\mathbb{P}_{t_j}} = \frac{\prod_{i=1}^n \exp\left\{-\frac{1}{2}\left(X_i-f_M(i)\right)^2\right\}}{\prod_{i=1}^n \exp\left\{-\frac{1}{2}\left(X_i-f_j(i)\right)^2\right\}} & =\exp\left\{\sum_{i=\floor{t_jn}+1}^n\frac{1}{2}\l(\left(X_i-f_j(i)\right)^2-\left(X_i-f_M(i)\right)^2\r)\right\},\\
      \text{and}\qquad  
 \frac{1}{2}\l(\left(X_i-f_j(i)\right)^2-\left(X_i-f_M(i)\right)^2\r)&=\frac{1}{2}\l(2X_i\left(f_M(i)-f_j(i)\right)+f_j^2(i)-f_M^2(i)\r)\\
        &=\varepsilon_i \l(f_M(i)-f_j(i)\r)-\frac{1}{2}\left(f_M(i)-f_j(i)\right)^2,
    \end{align*}
    where \(\varepsilon_i = X_i - f_j(i)\sim \mathcal{N}(0,1)\) if \(\left(X_i\right)_{i=1}^n\sim\mathbb{P}_{t_j}\).
    Now we can use the moment generating function of the standard normal distribution, \(\E\left[\exp\left\{tZ\right\}\right]=\exp\left\{t^2/2\right\}\),  to get
    \begin{align}
    &\E_{t_j}\left[\exp\left\{\varepsilon_i \l(f_M(i)-f_j(i)\r)-\frac{1}{2}\left(f_M(i)-f_j(i)\right)^2\right\}\right]= 1.\label{momentofone}
    \end{align}
    We set \(u_j \coloneqq t_j + b\psi_n\). Note that \(f_M(i)=0\) for \(i\le u_j < t_M\). We define \(g_n\) as
    \begin{equation}g_n\coloneqq\sum_{i=\floor{t_jn}+1}^{\floor{u_jn}} \frac{1}{2}\left(f_{M}(i)-f_{j}(i)\right)^2 = \frac{1}{2}\sum_{i=\floor{t_jn}+1}^{\floor{u_jn}} \left(\beta_+\l(\frac{i}{n}-t_j\r)+\alpha_+\right)^2= \frac{1}{2}\sum_{i=1}^{\floor{b\psi_nn}} \left(\beta_+\frac{i}{n}+\alpha_+\right)^2.\nonumber
    \end{equation}
    Note that \(g_n\) is not dependent on \(j\). Next we define \(Z_j\) as
    \[ Z_j \coloneqq \left(-\sum_{i=\floor{t_jn}+1}^{\floor{u_jn}}\varepsilon_if_j(i)\right)\left/{\sqrt{\sum_{i=\floor{t_jn}+1}^{\floor{u_jn}}f_j(i)^2}}\right.\sim \mathcal{N}(0,1).\]
    As \(\tilde{\tau}\) is a Markov stopping time,  \(\mathcal{B}_j= \left\{\left|\tilde{\tau}-t_j\right|\le b \psi_n\right\}\) is \(\mathcal{F}_{u_j}\)-measurable and is thus independent of \(X_i\)  and \(\varepsilon_i\), for \(i \in\{ \floor{u_jn}+1, \ldots, n\}\). 
   By this independence and \eqref{momentofone}, we rewrite the sum in \eqref{one} as
    \begin{align*}
        \mathbb{E}_{t_j}\left[\frac{\mathbf{d}\mathbb{P}_{t_M}}{\mathbf{d}\mathbb{P}_{t_j}} \mathbb{I}\left(\mathcal{B}_j\right)\right] &= \mathbb{E}_{t_j}\left[\exp\left\{\sum_{i=\floor{t_jn}+1}^n\varepsilon_i \l(f_M(i)-f_j(i)\r)-\frac{1}{2}\left(f_M(i)-f_j(i)\right)^2\right\} \mathbb{I}\left(\mathcal{B}_j\right)\right]\\
        &= \mathbb{E}_{t_j}\left[\exp\left\{\sum_{i=\floor{t_jn}+1}^{\floor{u_jn}}\varepsilon_i \l(f_M(i)-f_j(i)\r) -\sum_{i=\floor{t_jn}+1}^{\floor{u_jn}}\frac{1}{2}\left(f_M(i)-f_j(i)\right)^2\right\} \mathbb{I}\left(\mathcal{B}_j\right)\right]\\
        &\qquad\qquad\cdot\prod_{i = \floor{u_jn}+1}^n \E_{t_j}\left[\exp\left\{\varepsilon_i \l(f_M(i)-f_j(i)\r)-\frac{1}{2}\left(f_M(i)-f_j(i)\right)^2\right\}\right]\\
        &= \mathbb{E}_{t_j}\left[\exp\left\{\left(-\sum_{i=\floor{t_jn}+1}^{\floor{u_jn}}\varepsilon_i f_j(i)\right)-g_n\right\} \mathbb{I}\left(\mathcal{B}_j\right)\right]\\
        &= \mathbb{E}_{t_j}\left[\exp\left\{\sqrt{\sum_{i=\floor{t_j n}+1}^{\floor{u_j n}}\left(f_j(i)\right)^2} Z_j-g_n\right\} \mathbb{I}\left(\mathcal{B}_j\right)\right]\\
        &\ge \mathbb{E}_{t_j}\left[\exp\left\{\sqrt{\sum_{i=\floor{t_j n+1}}^{\floor{u_j n}}\left(f_j(i)\right)^2} Z_j-g_n\right\} \mathbb{I}\left(\mathcal{B}_j\right)\mathbb{I}(Z_j\ge 0)\right]\\
        &\ge \exp\left\{-g_n\right\}\mathbb{P}_{t_j}\left(\mathcal{B}_j\cap \left\{Z_j \ge 0\right\}\right).
    \end{align*}
    Because of \eqref{P(Bj)}, it holds that \(\mathbb{P}_{t_j}(\mathcal{B}_j) \ge 3/4\). Then
    \begin{align*}
        \mathbb{P}\left(\mathcal{B}_j \cap \left\{Z_j \ge 0\right\}\right) = \mathbb{P}(\mathcal{B}_j) + \mathbb{P}(Z_j\ge 0) - \mathbb{P}_{t_j}(\mathcal{B}_j \cup \left\{Z_j\ge 0\right\})
        \ge \frac{3}{4} + \frac{1}{2} - 1 = \frac{1}{4}.
    \end{align*}
    Combining this with \eqref{one}, we get
    \begin{equation}
        \frac{1}{4} \ge \frac{1}{4}\sum_{j=0}^{M-1} \exp\left\{-g_n  \right\} = \frac{M}{4}\exp\left\{-g_n\right\}
        = \floor{\frac{\l(1-2\delta_0\r)}{3b\psi_n}}\frac{1}{4} \exp\left\{ -g_n\right\}\asymp \frac{\exp\l\{-g_n\r\}}{\psi_n}.\label{four}
    \end{equation}
    \paragraph{For a jump}
    For a jump, we can set \(\beta_+=0\), \(|\alpha_+| \ge\delta_0\) and \(\psi_n = {\log n}/{n}\). 
    Then \[g_n = \frac{1}{2}\sum_{i=1}^{b\log n}\alpha_+^2 = \frac{\alpha_+^2b}{2}\log n. \]
    With \(b = \alpha_+^{-2}\), it holds that
    \begin{align*}
        \frac{n}{\log n}\exp\left\{-\frac{\alpha_+^2b}{2}\log n\right\}&=  \frac{n}{\log n}\exp\left\{-\frac{1}{2} \log n\right\}= \frac{\sqrt{n}}{\log n} \to \infty\text{, for } n \to \infty.
    \end{align*}
    Thus, \(\exp\left\{-g_n\right\}n/\log n\to\infty\) for \(n\to \infty\), which is a contradiction to \eqref{four}, so the statement is true.
    \paragraph{For a kink} Now we set \(\alpha_+ = 0\), \(|\beta_+|\ge\delta_0\) and  \(\left.\psi_n ={n^{1/3}}\middle/{\log^{{1}/{3}}n}\right. \). Then, it holds 
    \begin{align*}
        g_n &= \frac{1}{2}\sum_{i=1}^{bn^{2/3}\log^{1/3}n}\left(\beta_+\frac{i}{n}\right)^2 =\frac{\beta_+^2}{12n^2}(2b^3n^2\log n+3b^2n^{4/3}\log^{2/3}n+bn^{2/3}\log^{1/3}n)\\
        &= \frac{\beta_+^2b^3 \log n }{6}+\frac{\beta_+^2b^2\log^{{2}/{3}}n}{4n^{2/3}}+\frac{\beta_+^2b\log^{{1}/{3}}n}{12n^{4/3}}.
    \end{align*}
    For \(b = \beta_+^{-{2}/{3}}\), we get
    \[\exp\left\{-\frac{\beta_+^2b^2\log^{2/3}n}{4n^{2/3}}-\frac{\beta_+^2b\log^{{1}/{3}}n}{12n^{4/3}}\right\}\to 1 \text{, for } n\to\infty,\]
    and
    \[\frac{n^{1/3}}{\log^{1/3}n}\exp\left\{- \frac{\beta_+^2b^3\log n }{6}\right\} = \frac{n^{1/3}}{\log^{1/3}n} \exp\left\{-\frac{1}{6}\log n\right\} = \frac{n^{1/6}}{\log^{1/3}n}\to \infty \text{, for } n\to \infty.\]
    Thus, \(\bigl({n^{1/3}}/{\log^{{1}/{3}}(n)}\bigr)\exp\left\{-g_n\right\}\to\infty\), for \(n\to \infty\), which is a contradiction to \eqref{four}. 
\end{proof}

\section{Discussion}\label{s:dis}

We have studied the online detection of changes in segmented linear models with additive i.i.d.\ Gaussian noise. Our focus is on the minimax rate optimality in estimating the change point as well as computational and memory efficiency. We introduce the online detector FLOC, which offers several practical advantages, including ease of implementation as well as constant computational and memory complexity for every newly incoming data point --- crucial attributes for effective online algorithms. From a statistical perspective, FLOC achieves minimax optimal rates for detecting changes in both function values (i.e.\ jump) and slopes (i.e.\ kink). We believe that this is of particular practical benefit, as in many applications the type of change is not always clear beforehand. Notably, our results reveal a phase transition between the jump and kink scenarios, which echo the understanding in the offline setup (\citealp{GoTsZe06}, \citealp{FrHoMu14}, \citealp{Che21}; see also \Cref{tab:rates}). {\rev Furthermore, once a change point is detected, FLOC can, with high probability, distinguish between these two types of structural changes (a jump or a kink).} The FLOC detector is specifically designed to achieve asymptotically minimax optimal rates. While the constants involved have not been fully optimized and could likely be improved, we provide preliminary guidance for tuning FLOC to improve its empirical performance in finite-sample settings. Alternative approaches for parameter tuning could further enhance the performance of FLOC. For instance, theoretical insights, such as  the limiting distribution of detection delay provided by \citet{AuHoRe09}, could guide the selection of thresholds to satisfy specified bounds on type~II error, and bootstrapping methods, as introduced by \citet{HuKi12}, could be adapted to improve performance, particularly in small-sample scenarios. The current implementation of FLOC relies on sufficient historical data to accurately estimate the pre-change signal. As a practical extension, an adaptive approach could be developed to incrementally update the signal estimate as new observations become available. 

Monitoring simultaneously jumps and kinks can enhance detection power compared to conventional approaches that focus solely on mean changes, as  demonstrated in our analysis of excess mortality data. However, practitioners should be aware that the Gaussian noise and linear signal assumptions may be strongly violated in certain real-world applications.
{\rev Empirically, FLOC exhibits robustness under heavy tails.}
Nevertheless, enhancing the robustness of FLOC to accommodate broader noise distributions and signal structures, represents a promising direction for future research. For example, in the application discussed in \Cref{section:realdata}, it is  interesting to develop a refined analysis by incorporating serial dependencies. Finally, from a theoretical perspective, the extension of minimax rate results to more general settings, such as piecewise polynomial signals or changes in higher-order derivatives, remains an open problem and is worth further investigation. 


\section*{Acknowledgements}
The authors gratefully acknowledge support from the Deutsche Forschungsgemeinschaft (DFG, German Research Foundation) under Germany’s Excellence Strategy--EXC 2067/1-390729940. 
AM additionally acknowledges support from the DFG Research Unit 5381. 


%

\end{document}